\documentclass[11pt]{amsart}

\usepackage[all]{xypic}

\usepackage{latexsym}

\usepackage{amssymb}
\usepackage{amsfonts}
\usepackage{amscd}
\usepackage{amsmath,amsthm}

\newtheorem{lemma}{Lemma}[section]

\newtheorem{lem}[lemma]{Lemma}
\newtheorem{prop}[lemma]{Proposition}
\newtheorem{thm}[lemma]{Theorem}
\newtheorem{cor}[lemma]{Corollary}

{

}

\theoremstyle{definition}

\theoremstyle{remark}

\numberwithin{equation}{section}

\newenvironment{pf}{\noindent{\bf Proof.}}{\hfill $\square$\medskip}

\def\CC{{\mathbb C}}

\def\GG{{\mathbb G}}
\def\HH{{\mathbb H}}

\def\NN{{\mathbb N}}

\def\PP{{\mathbb P}}

\def\RR{{\mathbb R}}

\def\ZZ{{\mathbb Z}}

\def\bol{{\bar b}}

\def\iol{{\bar i}}
\def\jol{{\bar j}}
\def\kol{{\bar k}}

\def\ool{{\bar o}}

\def\uol{{\bar u}}
\def\vol{{\bar v}}

\def\0ol{{\bar 0}}
\def\1ol{{\bar 1}}
\def\2ol{{\bar 2}}
\def\ol2{{\bar 2}}
\def\3ol{{\bar 3}}
\def\4ol{{\bar 4}}
\def\5ol{{\bar 5}}
\def\6ol{{\bar 6}}
\def\7ol{{\bar 7}}
\def\8ol{{\bar 8}}
\def\9ol{{\bar 9}}

\def\bold0{{\bf 0}}
\def\bold1{{\bf 1}}
\def\bold2{{\bf 2}} 
\def\bold3{{\bf  3}}
\def\bold4{{\bf 4}}
\def\bold5{{\bf 5}}
\def\bold6{{\bf 6}}
\def\bold7{{\bf 7}}
\def\bold8{{\bf 8}}
\def\bold9{{\bf 9}}

\def\xul{{\underline{x}}}

\def\P2Skly{\PP^2_{Skly}}

\def\coker{\operatorname {coker}}

\def\End{\operatorname {End}}
\def\Ext{\operatorname {Ext}}

\def\GL{\operatorname {GL}}
\def\gr{\operatorname {gr}}

\def\Hom{\operatorname {Hom}}

\def\im{\operatorname {im}}

\def\ker{\operatorname {ker}}

\def\pd{{\operatorname {\partial}}}

\def\SU{\operatorname {SU}}

\def\th{\operatorname {th}}    
\def\Tor{\operatorname {Tor}}

\def\Tr{\operatorname {Tr}}

\def\Aut{\operatorname{Aut}}

\def\Bimod{{\sf Bimod}}

\def\fchar{\operatorname{char}}

\def\Der{\operatorname{Der}}

\def\dim{\operatorname{dim}}

\def\End{\operatorname{End}}

\def\Ext{\operatorname{Ext}}

\def\fdim{{\sf fdim}}

\def\gldim{\operatorname{gldim}}

\def\gr{{\sf gr}}

\def\Hom{\operatorname{Hom}}

\def\id{\operatorname{id}}

\def\Im{\operatorname{Im}}

\def\mod{{\sf mod}}

\def\pdim{\operatorname{pdim}}

\def\pd{{\partial}}

\def\Proj{\operatorname{Proj}}
\def\Projnc{\operatorname{Proj}_{nc}}

\def\proj{\operatorname{proj}}

\def\rank{\operatorname{rank}}

\def\Real{\operatorname{Re}}

\def\rep{{\sf rep}}

\def\RHom{\operatorname {RHom}}

\def\ul1{\operatorname{\underline{1}}}

\def\l{\leftarrow}

\def\d{\downarrow}

\def\a{\alpha}
\def\b{\beta}
\def\c{\gamma}
\def\d{\delta}

\def\l{\lambda}

\def\s{\sigma}

\def\ve{\varepsilon}

\def\L{\Lambda}

\def\fa{{\mathfrak a}}
\def\fb{{\mathfrak b}}
\def\fc{{\mathfrak c}}
\def\fcol{\overline{\mathfrak c}}

\def\fg{{\mathfrak g}}
\def\fgol{\overline{\mathfrak g}}

\def\fh{{\mathfrak h}}

\def\fsl{{\mathfrak s}{\mathfrak l}}
\def\fso{{\mathfrak s}{\mathfrak o}}

\def\sD{{\sf D}}

\def\wtA{{\widetilde{ A}}}

\def\cal{\mathcal}

\def\cA{{\cal A}}

\def\cC{{\cal C}}
\def\cD{{\cal D}}

\def\cL{{\cal L}}
\def\cM{{\cal M}}

\def\cO{{\cal O}}

\def\cT{{\cal T}}
\def\cU{{\cal U}}

\def\coh{{\sf coh}}

\def\dirlim{\mathop{\vtop{\baselineskip -100pt\lineskip -1pt\lineskiplimit 0pt
\setbox0\hbox{lim}\copy0\hbox to \wd0{\rightarrowfill}}}\limits}
\def\invlim{\mathop{\vtop{\baselineskip -100pt\lineskip -1pt\lineskiplimit 0pt
\setbox0\hbox{lim}\copy0\hbox to \wd0{\leftarrowfill}}}\limits}

\def\I11{{1 \kern -0.8pt \! \mbox{l}}}
\def\mumu{{\mu\kern-4.2pt\mu}}
\def\bfmu{{\mu\kern-4.2pt\mu}}
\def\2slash{\backslash \! \backslash}

\def\boxtimes{\setbox0\hbox{$\Box$}\copy0\kern-\wd0\hbox{$\times$}}

\def\udot{{\:\raisebox{2pt}{\text{\circle*{1.5}}}}}
\def\hdot{{\:\raisebox{2pt}{\text{\circle*{1.5}}}}}
\def\ldot{{\:\raisebox{2pt}{\text{\circle*{1.5}}}}}

\date{}                                           

\begin{document}

\title[A 3-Calabi-Yau algebra with $G_2$ symmetry]{A 3-Calabi-Yau algebra with $G_2$ symmetry
constructed from the octonions}

\author{S. Paul Smith}

\address{ Department of Mathematics, Box 354350, Univ.
Washington, Seattle, WA 98195}

\email{smith@math.washington.edu}

\thanks{S. P. Smith was supported by NSF grant DMS  0602347.}


\begin{abstract}
This paper concerns an associative $\NN$-graded algebra $A$ that is the enveloping algebra
of a Lie algebra with exponential growth. The algebra $A$ is 3-Calabi-Yau. There is a $\ZZ$-form of $A$ so for every field $k$ there is an algebra $A_k$. An algebraic group of type 
$G_2$ acts as degree-preserving automorphisms of $A$. 
We can define $A$ in many ways.

If  $V$ is the 7-dimensional irreducible representation of the complex semisimple Lie algebra of 
type $G_2$, then  $A_\CC$ is isomorphic to the tensor algebra $T(V)$ modulo the ideal generated by the submodule of $V^{\otimes 2}$ isomorphic to $V$. 

Alternatively, $A$ can be defined as a superpotential algebra derived from a 3-form on $\RR^7$ having an open $\GL(7)$ orbit and compact isotropy group. 
Or $A$ can be defined in terms of the product on the octonions.

Classification of the finite-dimensional representations of $A$ is equivalent to classifying square
matrices $Y$ with purely imaginary octonion entries such that $\Im(Y^2)=0$.

$A$ is a Koszul algebra. Its Koszul dual is a quotient of the exterior algebra on 7 variables, 
has Hilbert series $1+7t+7t^2+t^3$, and is a symmetric Frobenius algebra. 

Although not noetherian, $A$ is graded coherent, meaning that
the finitely presented graded $A$-modules form an abelian category. 
Let $\coh X$ be the quotient of the category of finitely presented graded $A$-modules by 
the full subcategory of finite dimensional graded modules.
The category $\coh X$ behaves somewhat like the category of coherent sheaves on the projective plane $\PP^2$. 
For example, analogous to a result for $\PP^2$, there are equivalences of bounded derived categories
$$
\sD^b(\coh X) \equiv \sD^b(\rep Q) \equiv \sD^b(\mod E) 
$$
where $Q$ is a quiver on three vertices (with relations--see section \ref{sect.nag})
having path algebra 
$$
E:= \begin{pmatrix}
  k  &A_1  & A_2   \\
   0 &  k  &A_1  \\   
   0 & 0 & k
\end{pmatrix}.
$$
An algebraic group of type $G_2$ acts as  automorphisms  of $E$. When $k$ is algebraically closed the indecomposable representations of dimension $(1,1,1)$ are parametrized by the tautological 
$\PP^2$-bundle over a smooth quadric hypersurface in $\PP^6=\PP(\Im {\mathbb O}_k)$. 
The classification of representations of $Q$ is equivalent to classifying pairs $(X,Y)$ of 
matrices with purely imaginary octonion entries such that $\Im(YX)=0$.

\end{abstract}

\maketitle

\pagenumbering{arabic}


\setcounter{section}{0}
 
 \vfill
 \eject
 
 \section{Introduction}

This paper initiates the study of an infinite dimensional associative algebra, $A$, that is intimately related to the exceptional  Lie group of type $G_2$. The algebra can be defined over any field but in this introduction we  focus on the real and complex versions, $A_\RR$ and $A_\CC=A_\RR \otimes_\RR \CC$.  

The structure constants for a suitable generating set  for $A_\RR$ belong to $\{0,\pm 1\}$ so there is a $\ZZ$-form of $A$ and hence an algebra $A_k$ defined over any field $k$. Some of the properties of $A_\RR$ and $A_\CC$ stated in this introduction hold for all $A_k$ and in that case we will simply write $A$. For example, $A$  is the enveloping algebra of an infinite dimensional Lie algebra. 

\subsection{The origin of $A$}
The algebra $A$ was brought to my attention by Richard Eager who asked if it is Koszul. It is. Eager had 
observed that a superpotential in the paper \cite{FST} would, through taking cyclic partial derivatives, give rise to an algebra $A$. The superpotential in \cite{FST} is constructed in a straightforward way from a generic 3-form on $\CC^7$. The algebra $A$ is not mentioned in \cite{FST}.

 \subsection{Definition of $A_\RR$ in terms of the octonions}
 \label{sect.defn1}
 One of the simplest ways to define $A_\RR$ is the following. 
Consider the diagram 
$$
  \UseComputerModernTips
\xymatrix{
{\mathbb O} \otimes {\mathbb O} \ar[rr]^{\rm mult} && {\mathbb O} \ar[d]^{\Im}
\\
\Im{\mathbb O} \otimes \Im{\mathbb O} \ar[u]^{\iota \otimes \iota} \ar[rr]_{\mu} && \Im{\mathbb O} 
}
$$
where the top arrow is multiplication, $\iota$ is the inclusion of the purely imaginary octonions,
and $\mu$ the composition of the other three arrows. Let
$$
\mu^* : (\Im{\mathbb O})^* \longrightarrow  (\Im{\mathbb O})^* \otimes (\Im{\mathbb O})^*
$$
be the dual of $\mu$ and define
$$
A_\RR:={{T \Im{\mathbb O}^*}\over{(\im\mu^*)}}.
$$
Here $T  \Im{\mathbb O}^*$ denotes the tensor algebra on the vector space dual to $\Im {\mathbb O}$ and 
$(\im \mu^*)$ is the two-sided ideal generated by the image of $\mu^*$. Thus $A$ is generated by seven elements modulo seven relations. 

The action of the compact simple real Lie group $G_2^c$ as automorphisms of ${\mathbb O}$ induces an action of $G_2^c$ as algebra automorphisms of $A_\RR$. 

If $u,v \in {\mathbb O}$, then $\Im(uv)=-\Im(vu)$ so the image of $\mu^*$ consists of skew-symmetric tensors.
Thus $A_\RR$ maps onto the polynomial ring $S(\Im{\mathbb O}^*)$.

\subsection{$A_\RR$ as an enveloping algebra}
\label{sect.env.alg}
Let ${\mathfrak f}$ be the free Lie algebra on the vector space $\Im{\mathbb O}^*$. Then $T  \Im{\mathbb O}^*$ is 
the enveloping algebra of $\mathfrak f$. Because $\im \mu^*$ consists of skew-symmetric tensors 
it is contained in ${\mathfrak f}$. Let $\overline{\mathfrak f}$ be the quotient of ${\mathfrak f}$ by the Lie ideal generated by 
$\im \mu^*$. Then $A$ is the enveloping algebra of  the Lie algebra $\overline{\mathfrak f}$.    

\subsection{Definition of $A_\CC$ in terms of Lie theory}
Let $\fg_2$ be the simple complex Lie algebra of type $G_2$. Let $V$ be the unique 7-dimensional 
irreducible representation of $\fg_2$ and let $R$ be the  unique subspace of $V \otimes V$ that is isomorphic to $V$ as a $\fg_2$-module. Then
$$
A_\CC \cong {{TV}\over{(R)}}.
$$
It is well-known that $R$ is contained in the space of skew-symmetric tensors so, independently of the observation that $\Im(uv)=-\Im(vu)$ in section \ref{sect.defn1}, the ``same'' argument as in section \ref{sect.env.alg} shows that $A$ is the enveloping algebra of a Lie algebra. 

\subsection{Construction of $A$ from a generic 3-form on $\RR^7$}
\label{sect.G2.holonomy}
There are two open orbits for the action of $\GL(7)$ on $\wedge^3 \RR^7$. The isotropy groups of the points
in one of those orbits are isomorphic to the compact simple Lie group of type $G_2$. Let $\phi$ belong to that orbit. The 3-form $\phi$ plays a fundamental role in the theory of manifolds with $G_2$-holonomy. 
Let $W$ be the unique totally anti-symmetric element of $V^{\otimes 3}$
mapping onto $\phi$. Then $A_\RR$ is isomorphic to the super-potential algebra associated to $W$, i.e., $A_\RR$ is the tensor algebra on $\RR^7$ modulo the relations given by the cyclic partial derivatives of $W$. 
(Section \ref{sect.spp.alg} and Proposition \ref{prop.spp.alg}).

\subsection{Properties of $A$}

The algebra $A$ is $\NN$-graded with degree-one component equal to $\Im{\mathbb O}^*$. Its Hilbert series,  $\sum (\dim A_n)t^n$,  is $(1-7t+7t^2-t^3)^{-1}$. In particular, $A$ has exponential growth---it grows faster than the free algebra on five variables but slower than the free algebra on six variables. It comes as no surprise that 
$A$ is not noetherian. It is the enveloping algebra of an infinite dimensional Lie algebra of exponential growth. 

However, 
$A$ is a coherent algebra: its category of finitely presented modules in abelian. This was first proved by Dmitri Piontkovskii and I am grateful for his allowing me to include his proof here.

 In section \ref{sect.4.2} we show that $A$'s homological properties are as good as those of the commutative polynomial ring in 3 variables. 
The global homological dimension of $A$ is 3. In fact, $A$ is a Koszul algebra and the homogeneous components of the Koszul dual $\Ext^\hdot_A(\RR,\RR)$ have dimensions 1,7,7,1. Because the polynomial ring on seven variables is a quotient of $A$, $\Ext^\hdot_A(\RR,\RR)$ is a quotient of the exterior algebra
$\Lambda^\hdot \RR^7$. 

More importantly, $A$ is a 3-dimensional Calabi-Yau algebra; i.e.,  the natural map 
$A \to \Ext^3_{\sf Bimod}(A, A \otimes A)$ is an isomorphism of $A$-bimodules and $\Ext^i_{\sf Bimod}(A, A \otimes A)=0$ when $i \ne 3$. It is striking that $A_\RR$ is a 3-dimensional Calabi-Yau algebra and also 
can be obtained from the 3-form that appears in the theory of manifolds with $G_2$-holonomy. This combination of $G_2$ and Calabi-Yau  echos the fact that string theory requires a 10-dimensional space fibered  over $\RR^{3,1}$ by compact Calabi-Yau 3-folds whereas M-theory requires an 
11-dimensional space fibered  over  $\RR^{3,1}$ by compact manifolds with $G_2$-holonomy.

\subsection{Representations of $A$}
Although we don't know much about the representation theory of $A_\RR$ the finite dimensional representation theory of $A$ is equivalent to the following intriguing, and deceptively innocent-sounding, problem: classify, up to conjugation by 
$\GL(n,\RR)$,  all $n \times n$ matrices $X$ with entries in ${\mathbb O}$ such that the purely imaginary part of $X^2$ is zero.

Because the quotient of $A$ modulo the ideal generated by all commutators is the symmetric algebra $S(\Im{\mathbb O}^*)$, the 1-dimensional representations of $A_\RR$ are parametrized by $\Im {\mathbb O}$. It is somewhat surprising that if $M$ and $N$ are non-isomorphic 1-dimensional representations of $A_\RR$, then $\Ext^1_A(M,N)=0$. The quaternion algebra is not a quotient of $A_\RR$. The   $n \times n$ matrix ring $M_n(\CC)$ is a quotient of  $A_\RR$ for all $n \ge 1$. The real Lie algebra $\fso(3,1)$ is a quotient of $\overline{\mathfrak f}$ so every representation of $\fso(3,1)$ is an $A$-module. For example, $M_4(\RR)$ is a quotient of $A$. 
All 2-dimensional irreducible representations of $A_\RR$ arise from homomorphisms $A_\RR \to \CC$ so 
$M_2(\RR)$ is not a quotient of $A$. We do not know if $M_3(\RR)$ is a quotient of $A$.

\subsection{Quotients and subalgebras of $A$}
The free algebra on two variables, $\CC\langle x,y\rangle$, is a quotient of $A_\CC$, but $\RR\langle x,y\rangle$ is not a quotient of $A_\RR$. The quadratic self-duality algebra of Connes and Dubois-Violette 
\cite{CDV}, $\cA^{(+)}$ in their notation, is a quotient of $A$. There is a striking dichotomy involving 
$\cA^{(+)}$
which is best stated after identifying $(\Im {\mathbb O})^*$ with $\Im {\mathbb O}$ via the symmetric bilinear 
form $\langle u,v\rangle =\Real(u\vol)$ on ${\mathbb O}$; this identification allows us to realize $A$ as a 
quotient of $T(\Im {\mathbb O})$ and identify the degree-one component of $A$ with $\Im {\mathbb O}$. The dichotomy is this: if $L$ is a 3-plane in $\Im {\mathbb O}$, then $A/(L)$ is isomorphic to $\cA^{(+)}$
if $L = \Im \HH$ for some copy of the quaternions $\HH \subset {\mathbb O}$ and is isomorphic to the 
polynomial ring on 4 variables otherwise. Subspaces of the form $\Im \HH$ are called associative 3-planes 
in the differential geometry literature and play an important role in a wide range of matters. A subspace of $\Im {\mathbb O}$ orthogonal to an associative 3-plane is said to be co-associative. The key step in Piontkovskii's proof that $A$ is coherent is that the 2-sided ideal generated by a co-associative 4-plane is a free 
$A$-module and the quotient by that ideal is a polynomial ring in three variables. 

Let $B$ be a subalgebra of $A_\RR$ generated by any 6-dimensional subspace of the degree-one component of $A_\RR$. Then algebra $A_\RR$ is an Ore extension of $B$ with respect to a suitable derivation. All such $B$ are isomorphic because $G_2^c$ acts transitively on the lines in $\Im
{\mathbb O}$ and hence on their orthogonals because the bilinear form $\langle -,-\rangle$ is
$G_2^c$-invariant.
The algebra $B$ is a free algebra modulo one homogeneous quadratic relation of rank equal to $\dim B_1$. General results of James Zhang \cite{Z} imply that $\gldim B=2$ and $B$ is regular in the sense that $\Ext^i_B(\RR,B)$ is zero when $i \ne 2$ and is $\RR$ when $i=2$. We  show that 
$A_\RR$ has no  degree-two relations of rank $<6$ and then use results in \cite{Z} to show that the center of $A$ is equal to $\RR$. In fact, $A_\RR$ has no normal elements except those in $\RR$. 

The algebra $B$ has a close relationship to the preprojective algebra for the generalized Kronecker quiver with 6 arrows.

\subsection{Non-commutative geometry and the derived category}

Let $\coh A$ denote the abelian category of finitely presented graded $A$-modules modulo the subcategory
of finite dimensional modules. We implicitly define a non-commutative variety $X$ by declaring that the ``category of coherent sheaves'' on $X$, which we denote by $\coh X$, is the quotient of $\coh A$ by the localizing subcategory consisting of finite dimensional graded modules. There are equivalences of bounded derived categories
$$
\sD^b(\coh X) \equiv \sD^b(\rep Q) \equiv \sD^b(\mod R) 
$$
where $Q$ is a quiver
$$
 \UseComputerModernTips
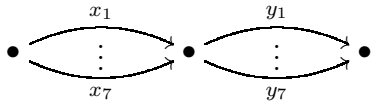
\xymatrix{
\bullet \ar@/^.8pc/[rr]^{x_1}  \ar@/_.8pc/[rr]_{x_7}^{\vdots} && \bullet \ar@/^.8pc/[rr]^{y_1}   \ar@/_.8pc/[rr]^{\vdots}_{y_7}  && \bullet
}
$$
with 3 vertices, 14 arrows, and 7 relations, and  path algebra
$$
E= \begin{pmatrix}
  \RR  & (\Im {\mathbb O})^*  & \coker \mu^*   \\
   0 &  \RR  &  (\Im {\mathbb O})^*   \\   
   0 & 0 & \RR
\end{pmatrix}. 
$$
The compact Lie group $G^c_2$ acts as automorphisms of $E$. 

Write $\PP^6$ for the projective space of lines in $\Im {\mathbb O}_k$. 
We show that the moduli space of indecomposable representations of $Q$ of dimension $(1,1,1)$ is 
the subvariety of $\PP^6 \times \PP^6$ consisting of points $([u],[v])$ such that $\Im(uv)=0$. Over $\CC$
this has two components, the diagonal copy of $\PP^6$ and a $\PP^2$-bundle over a smooth quadric hypersurface in $\PP^6$.

 \subsection{Acknowledgments}
I am particularly grateful to Richard Eager for directing my attention to this algebra  after seeing
equations (14) and (15) in \cite{FST}. Dmitri Piontkovski proved Proposition 
\ref{prop.coherent} and kindly allowed me to include it here. Izuru Mori told me that the equivalence of categories in  Theorem \ref{thm.equiv} follows from his work with Minamoto. Monty McGovern, Michel Van den Bergh, and James Zhang, made helpful remarks. I thank them all. Part of this work was done during a visit to the Universidad de Buenos Aires and I am grateful to Andrea Solotar for her hospitality and conversations about this work.

\section{Definitions of $A$}

There are many ways to define $A$ and the different definitions emphasize and suggest different features of $A$. Ultimately, though, the variety among the different definitions rests on the fact that the Lie group of type $G_2$ can be defined in many ways. 
 
 To show the equivalence of the various definitions we need to start with an ``official'' definition and then
 show the other definitions produce an isomorphic algebra. Our official definition of $A$ over a field $k$ will be given in terms of ``the'' octonion algebra over $k$. Thus, the official definition of $A_\RR$ is that given in section \ref{sect.defn1} above. In section \ref{sect.oct} we will choose a basis for ${\mathbb O}$ with respect to which
 the structure constants belong to $\{0,\pm 1\}$.  The $\ZZ$-span of that basis is therefore a subalgebra of 
 ${\mathbb O}$. We denote it by ${\mathbb O}_\ZZ$, define ${\mathbb O}_k := {\mathbb O}_\ZZ 
 \otimes_\ZZ k$, and  then define $A_k$ as we defined $A_\RR$ in section \ref{sect.defn1}. That will be the official definition of $A_k$. Alternatively,  the $\ZZ$-form of ${\mathbb O}$ leads to a $\ZZ$-form of $A$, namely $A_\ZZ$, and we can define $A_k = A_\ZZ \otimes_\ZZ k$. 
 
$A_k$ can also be defined in terms of generators and relations (Proposition \ref{prop.relns1}).
 Proposition \ref{prop.relns/C}  gives a ``better'' set of generators and relations for $A_\CC$ that is closely 
 related to the Pfaffian system that Cartan used in his ``five variables paper''  \cite{EC1} to characterize the complex Lie group of type $G_2$.

\subsection{The official definition of $A_k$}
\label{sect.gens.relns}
As in section \ref{sect.defn1}, we define $A_k$ to be the tensor algebra on the dual $(\Im {\mathbb O}_k)^*$
modulo the ideal generated by the image of the map $\mu^*$ where $\mu$ is the map in the commutative diagram 
\begin{equation}
\label{diag.defn.A}
  \UseComputerModernTips
\xymatrix{
{\mathbb O}_k \otimes {\mathbb O}_k \ar[rr]^{\rm mult} && {\mathbb O}_k \ar[d]^{\Im}
\\
\Im{\mathbb O}_k \otimes \Im{\mathbb O}_k \ar[u]^{\iota \otimes \iota} \ar[rr]_{\mu} && \Im{\mathbb O}_k
}
\end{equation}
We make $A_k$ a graded $k$-algebra by placing $\Im{\mathbb O}^*$ in degree one.

\begin{prop}
There is an injective group homomorphism 
\begin{equation}
\label{Aut.A}
\Aut {\mathbb O}_k  \to \Aut_{\sf gr} A_k
\end{equation}
to the group of degree-preserving $k$-algebra automorphisms of $A_k$. 
\end{prop}
\begin{pf}
The map $\mu: \Im {\mathbb O}_k \otimes \Im {\mathbb O}_k \to \Im {\mathbb O}_k$ in 
(\ref{diag.defn.A}) is $\Aut {\mathbb O}_k$-equivariant because the fact that the decomposition ${\mathbb O}_k = k.1 \oplus \Im {\mathbb O}_k$ is a decomposition into $\Aut {\mathbb O}_k$-modules 
(section \ref{sect.oct}) implies that
 the three clockwise-directed arrows in (\ref{diag.defn.A}) are $\Aut {\mathbb O}_k$-equivariant.
So $\mu^*$ is $\Aut {\mathbb O}_k$-equivariant.
 
Every automorphism of $  {\mathbb O}_k$ therefore induces a $k$-linear automorphism of $(\Im {\mathbb O}_k)^*$ that extends in a unique way to a $k$-algebra automorphism of $A_k$. The induced automorphism preserves  each homogeneous component of $A_k$.

The only algebra automorphism of  
${\mathbb O}_k$ that restricts to the identity on $\Im {\mathbb O}_k$ is the identity map so the map in (\ref{Aut.A}) is  injective.
\end{pf}

\begin{thm}
\cite[Thm. 2.3.5]{SV}.
Let $\kol$ be an algebraic closure of $k$. Let ${\bf G}=\Aut ({\mathbb O}_{\kol})$. Then ${\bf G}$ is a connected, simple algebraic group of type $G_2$, and $\Aut {\mathbb O}_k={\bf G}(k)$,
the group of $k$-rational points. 
\end{thm}

\subsection{First properties of $A_k$}
Let $\mathfrak f$ denote the free Lie algebra over $k$ on the vector space $(\Im{\mathbb O}_k)^*$. We impose a grading on $\mathfrak f$ by setting $\mathfrak f_1=\Im {\mathbb O}_k^*$. 
Let $\overline{\mathfrak f}$ be the quotient of $\mathfrak f$ by the Lie ideal generated by $\im \mu^*$. This makes sense: if $u,v \in \Im{\mathbb O}_k$, then $\Im(uv)=-\Im(vu)$ so $\im \mu^*$ consists of skew-symmetric tensors and is 
therefore contained in ${\mathfrak f}_2$. 

\begin{prop}
\label{prop.A1}
$\phantom{xxx}$
\begin{enumerate}
  \item 
$A_k$ is the enveloping algebra of the Lie algebra $\overline{\mathfrak f}$.
\item{}
$A_k$ is a domain, i.e., it has no zero-divisors.
\end{enumerate}
Hence $A_k$ is  a Hopf algebra and isomorphic to its opposite algebra.
\end{prop}
\begin{pf}
For any vector space $V$, the tensor algebra $TV$ is isomorphic to the enveloping algebra of the free Lie algebra on $V$. We apply this with $V= (\Im{\mathbb O}_k)^*$ and note that if $\fh$ is an ideal in a Lie algebra $\fg$, then $U(\fg/\fh) \cong U(\fg)/(\fh)$. This proves (1), and (2) follows from the Poincar\'e-Birkhoff-Witt Theorem. 
\end{pf}

 \subsection{The Fano plane}
 We will use the Fano plane to give an explicit multiplication for the octonions. The Fano plane also serves as a mnemonic device for the relations for $A$ in Proposition \ref{prop.relns1}. 
 The Fano plane also helps us organize other properties of $A$.
 
The   {\sf Fano plane}  is the projective plane over the field of two elements.    
We label its points by the integers $1,\ldots,7$ in such a way that its seven 
lines consist of the triplets of points
\begin{equation}
\label{F.lines}
123,   \quad 145,  \quad 167,  \quad 246, \quad 275,  \quad  374,  \quad 365.
\end{equation}
We impose a cyclic ordering on each line by drawing an arrow from each point on the line to the ``next'' 
point on it. Below is the standard diagram of the Fano plane with most of  the arrows: we omit the arrows $1 \to 6$, $2 \to 3$, $3 \to 6$, $4 \to 3$, $5 \to 2$, and $6 \to 2$,
to avoid clutter:
 \begin{equation}
 \label{FP}
  \UseComputerModernTips
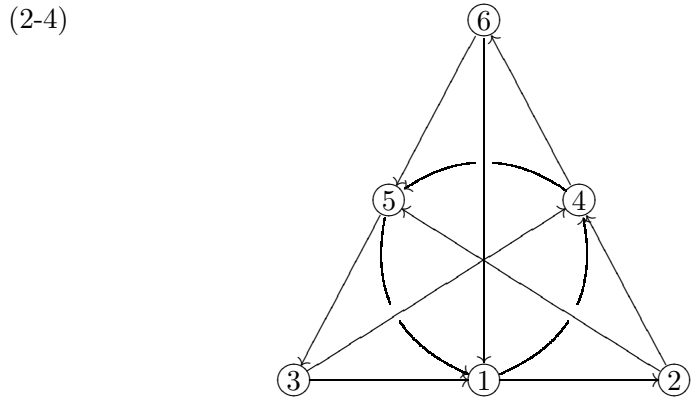
\xymatrix{
&& *+[o][F]{6} \ar[ddl]   \ar[dddd]    && \\
\\
& *+[o][F]{5} \ar[ddl]  \ar@/_1.5pc/[ddr] | \hole && *+[o][F]{4} \ar[uul]  \ar@/_1.2pc/[ll] | \hole &
 \\
&&   
  &&\\
*+[o][F]{3}  \ar[rr]   \ar[uurrr]   && *+[o][F]{1}  \ar@/_1.5pc/[uur] | \hole \ar[rr] & & *+[o][F]{2} \ar[uul]  \ar[uulll]
}
\end{equation}
The centroid is labelled 7. We have not drawn the heads of the arrows from the points 2, 6, and 3, that terminate at the point 7. 
If there are arrows $i \to j \to k \to i$ we call $ijk$ a {\sf directed line}. Thus (\ref{F.lines}) is a list of the directed lines. 

Mariano Su\'arez-Alvarez pointed out  that the organization of the arrows is the unique one up 
to automorphisms such that  the endpoints of the three arrows emanating from each point 
are colinear. For example, the endpoints of the arrows emanating from 6 are the colinear points 5,7,2;
the endpoints of the arrows emanating from 7 are 1,4,5.

\subsubsection{The symbol $\ve^{ijk}$.}
\label{ssect.ve}
If $\{i,j,k\} \subset \{1,2,3,4,5,6,7\}$ we define
$$
\ve^{ijk}:= 
\begin{cases}
+1 & \text{if $ijk$ is a directed line  in the Fano plane}
\\
\phantom{+} 0 & \text{if  $\vert\{i,j,k\}| \le 2$, or $i,j,k$ are distinct non-colinear}
\\
-1 & \text{if $jik$ is a directed line in the Fano plane.}
\end{cases}
$$
Explicitly,
$$
\begin{array}{l}
\ve^{123}=\ve^{145}=\ve^{167}= \ve^{246}=\ve^{275}=\ve^{374}=\ve^{365}=1,  \quad \;  \hbox{and }
\\
\ve^{ijk} = \ve^{kij} = - \ve^{jik}  \qquad \hbox{for all $i,j,k$}.
\end{array}
$$
Under the action of the symmetric group $S_3$ on $(i,j,k)$, $\ve^{ijk}$ transforms according to the sign representation.

\subsection{The octonions}
\label{sect.oct}

Let $k$ be a field of characteristic $\ne 2$. 

The {\sf octonion algebra} over $k$,  which we denote by ${\mathbb O}_k$, or ${\mathbb O}$, 
is the 8-dimensional, unital, non-associative, $k$-algebra with basis $\{1, o_1,\ldots,o_7\}$ and multiplication 
defined as follows: for each directed line $ijk$ in the Fano plane    $o_i^2=o_j^2=o_k^2=-1$ and 
  $$
  o_io_j = o_k = - o_jo_i \quad \hbox{and cyclic permutations}.
  $$
Equivalently, for any $r,s \in \{1,\ldots,7\}$, 
$$
o_ro_s=\sum_{i=1}^7 \ve^{rsi} o_i - \d_{rs}.
$$

If $ijk$ is a directed line, then the span of $1,o_i,o_j,o_k$ is a quaternion algebra and $\ve^{ijk}o_k=o_io_j$. 
Octonion algebras can be defined in greater generality \cite{SV} but only this definition is relevant to this paper. 
  
We write $\Im {\mathbb O}_k$ for the span of $o_1,\ldots,o_7$ and call it the space of {\sf imaginary octonions}.
{\sf Conjugation}, $u \mapsto \uol$, on  ${\mathbb O}_k$ is the $k$-linear map 
such that $\overline{1}=1$ and $\ool_i=-o_i$ for $i=1,\ldots,7$.  Thus $$\overline{uv} = \vol\uol.$$

Table  (\ref{mult.table.O}) is the multiplication table for the imaginary 
  octonions written with the following conventions: 
an entry $i$ in the table denotes $o_i$; $\iol$ denotes $-o_i$; a blank denotes $-1$ (for example $o_3o_3=-1$); the entry in row  $x$ and column $y$  is $xy$ (for example $o_5o_2=o_7$ and $o_2o_6=-o_4$).
  \begin{equation}
  \label{mult.table.O}
\begin{array}{|c|c|c|c|c|c|c|c|c|c|c|}
\hline
 & o_1 &  o_2  &  o_3 &  o_4  &  o_5 &  o_6  &  o_7   \\
\hline
  o_1 &   &    3     &  \ol2 & 5  & \4ol & 7 & \6ol \\
  o_2 &  \3ol &  & 1 &  6 & \7ol & \4ol & 5 \\
  o_3 & 2  & \1ol &  & \7ol  & \6ol & 5 & 4 \\
  o_4 &  \5ol & \6ol & 7 &   & 1 & 2 & \3ol \\
  o_5 &4  & 7 & 6 & \1ol  &  & \3ol &  \ol2 \\
 o_6  & \7ol & 4 & \5ol & \ol2  & 3 &  & 1 \\
  o_7 & 6 & \5ol & \4ol & 3  & 2 & \1ol &  \\
\hline
\end{array}
\end{equation}

Let $u \in {\mathbb O}_k$. The {\sf imaginary part} of $u$ is $\Im(u):={{1}\over{2}}(u-\uol)\in \Im {\mathbb O}_k$
 and   the {\sf real part} of $u$ is $\Real(u):={{1}\over{2}}(u+\uol) \in k.1$. 
We have
 $$
 \Im(u \vol)=-\Im(v\uol) \qquad \hbox{and} \qquad \Real(u\vol)=\Real(v\uol).
 $$
 If $u,v \in  \Im {\mathbb O}_k$, then 
 $\Im(uv)=-\Im(vu)$ and $\Real(uv)=\Real(vu)$.

\subsection{The non-degenerate symmetric bilinear form $\langle -,- \rangle$ on ${\mathbb O}_k$}

There is a non-degenerate symmetric bilinear form $\langle -,- \rangle$ on ${\mathbb O}_k$ defined by 
\begin{equation}
\label{eq.inn.prod}
\langle u,v \rangle:=\Real(u\vol).
\end{equation}
It is clear that $\{1,o_1,\ldots,o_7\}$ is an orthonormal basis for ${\mathbb O}_k$. In particular, 
$k.1$ and $\Im {\mathbb O}_k$ are  orthogonal to each other.
 
Let $a_0, \ldots,a_7 \in k$. Let $u=a_0+\sum_{i=1}^7 a_io_i$. Since $o_io_j=-o_jo_i$ for $i \ne j$,  
 \begin{align*}
u^2 & =a_0^2-\sum_{i=1}^7a_i^2+2a_0\sum_{i=1}^7a_io_i \qquad \hbox{and}
\\
u\uol &=a_0^2+\sum_{i=1}^7a_i^2 = \langle u,u\rangle.
\end{align*}
Therefore
\begin{equation}
\label{min.eq}
u^2 -2 \langle u,1 \rangle u  + \langle u,u \rangle =0.
\end{equation}
This identity implies that
\begin{enumerate}
  \item 
  $u \in \Im{\mathbb O}_k \qquad \qquad \Longleftrightarrow \qquad u^2 =-\langle u,u\rangle$;
  \item 
   $u \; \hbox{is a unit} \qquad \qquad  \Longleftrightarrow \qquad \langle u,u\rangle \ne 0$;
  \item 
     $u \; \hbox{is a zero-divisor} \phantom{...} \Longleftrightarrow \qquad \langle u,u\rangle = 0$.
\end{enumerate}
When ${\mathbb O}_k$ has zero-divisors we say it is {\sf split}. Note that 
(3) implies ${\mathbb O}_k$ is split if and only if  $-1$ is a sum of $\le 6$ squares in $k$. 

\begin{lem}
\label{g2.mods}
\label{lem.inv.form}
$\phantom{xxx}$
\begin{enumerate}
  \item 
  Both  $k.1$ and  $\Im{\mathbb O}_k$ are stable under the action of  $\Aut {\mathbb O}_k $.
  \item 
  The bilinear form $\langle -,-  \rangle$ is  $\Aut {\mathbb O}_k $-invariant.
  \item 
  If $\s \in \Aut {\mathbb O}_k $ and $u \in  {\mathbb O}_k $, then $\s\uol=\overline{\s u}$.
\end{enumerate}
\end{lem}
\begin{pf}
Let $u\in \Im {\mathbb O}_k$ and $\s \in \Aut {\mathbb O}_k $.

(1)
Applying (\ref{min.eq}) to $\s u$ in place of $u$ gives
$$
(\s u)^2 -2\langle \s u, 1\rangle \s u + \langle \s u,\s u \rangle =0.
$$
Applying $\s^{-1}$ to the previous equality gives
$$
u^2 -2\langle \s u, 1\rangle u + \langle \s u,\s u \rangle. =0.
$$
Comparing this equality with (\ref{min.eq}) gives 
$$
2(  \langle u,1 \rangle - \langle \s u,1 \rangle ) u = \langle u, u \rangle - \langle \s u, \s u \rangle.
$$
The right-hand side of this equality is in $k.1$ so if $u \notin k.1$ it follows that 
$ \langle\s u,1 \rangle =  \langle u,1 \rangle$ and 
$\langle u, u \rangle = \langle \s u, \s u \rangle$. In particular, $(k.1)^\perp$ is stable under the action of $\s$.
This proves (1)

(2)
The equality $\langle u, u \rangle = \langle \s u, \s u \rangle$ holds for all $u \in k.1$, so holds for all $u \in 
{\mathbb O}_k$.  Replacing $u$ by $u+v$ in this equality and expanding  gives 
$ \langle u, v \rangle = \langle \s u, \s v \rangle$.

(3)
Since $\s$ acts as the identity on $k.1$ and preserves $\Im {\mathbb O}_k$, 
$$
\s u + \overline{\s u}=2\Real(\s u) = 2 \Real u = 2\s(\Real u) = \s u + \s \uol
$$
so $\overline{\s u} =  \s \uol$.
\end{pf}

\subsection{The alternating trilinear form on $\Im {\mathbb O}_k$}
\label{sect.trilin.form}
Let $\{x_1,\ldots,x_7\}$ be the  basis for $(\Im {\mathbb O}_k)^*$ that is dual to $\{o_1,\ldots,o_7\}$ 
  with respect to the non-degenerate form $\langle -,-\rangle$. Thus  
  \begin{equation}
  \label{defn.xi}
  x_i:=\langle o_i,-\rangle 
  \end{equation}
With respect to the bilinear form $\langle -,-\rangle$ on $\Im {\mathbb O}_k$, $\{o_1,\ldots,o_7\}$ is an orthonormal basis.
  
Let $x^{ijk}=x_i \wedge x_j \wedge x_k \in \wedge^3 ( \Im {\mathbb O}_k)^*$.

\begin{prop}
\label{prop.3-form}
The formula
\begin{equation}
\label{eq.phi}
\phi(u,v,w):= -\Real(uvw)
\end{equation}
defines an $\Aut {\mathbb O}_k $-invariant alternating trilinear form on $\Im {\mathbb O}_k$. In terms of the 
dual basis,
\begin{align}
\label{phi}
\phi &= x^{123}+ x^{145}+x^{167} + x^{246}+x^{275}+x^{374}+x^{365} 
\\
& =
 \hbox{${{1}\over{6}}$} \sum_{i,j,k=1}^7 \ve^{ijk} x^{ijk}
 \end{align}
 where  $\ve^{ijk}$ is defined in section \ref{ssect.ve}.  
\end{prop}
\begin{pf}
The invariance follows from the fact that  $\Aut {\mathbb O}_k $ acts as the identity on $k$.  
Because $u$, $v$, and $w$, are purely imaginary, 
$$
\Real(uvw)=\Real(\Im(uv)w)=-\Real(\Im(vu)w) =- \Real(vuw).
$$
Likewise, $\Real(uvw) =-\Real(uwv)$, so $\phi$ is alternating.

Consider $\phi(o_i,o_j,o_k)$. If $|\{i,j,k\}| \ne 3$, then $\phi(o_i,o_j,o_k)=0$. Suppose  $|\{i,j,k\}| = 3$.
If $i,j,k$ are not colinear, then $o_io_jo_k \in \Im {\mathbb O}_k$ so $\phi(o_i,o_j,o_k)=0$.  
If $ijk$ is a directed line in the Fano plane, then $\phi(o_i,o_j,o_k)=1$. Hence $\phi$ is given by the formula 
(\ref{phi}).
\end{pf}

The form $\phi$ is the negative of the associative 3-form on \cite[p. 113]{H}.

The terms in $x^{123}+ x^{145}+x^{167} + x^{246}+x^{275}+x^{374}+x^{365}$ 
 are indexed by the directed lines in the Fano plane.

The subgroup of $\GL(7,\CC)$ that fixes $\phi_\CC$  is a simple Lie group of type $G_2$ \cite{En00},
  \cite[Defn. 1, p. 539]{B1}, \cite{Sch}.
The subgroup of $\GL(7,\RR)$ that fixes $\phi_\RR$ is isomorphic to $G_2^c$, the   compact, simple  real Lie group of type $G_2$ \cite{R}. Engel proved the result over $\CC$ in 1900 and his student, Reichel, proved the result over $\RR$ in 1907. Agricola gives a nice  account of their work \cite{A}.

  \subsection{Explicit relations for $A_k$}

  We define 
  \begin{equation}
  \label{defn.ri}
 r_i:=\mu^*(x_i)
  \end{equation}
  where $\mu^*$ is the dual to the map $\mu$ in (\ref{diag.defn.A}). Thus 
  $$
A_k = {{k \langle x_1,\ldots,x_7\rangle}\over{(r_1,\ldots,r_7)}}.
$$

\begin{prop}
\label{prop.relns1}
The algebra $A_k$ is $k\langle x_1,\ldots,x_7\rangle$ modulo the relations
$$
 \begin{array}{ cc}
 \begin{array}{c}
r_1=  [x_2,x_3]+[x_4,x_5]+[x_6,x_7]    \phantom{.} \\
r_2=    [x_3,x_1]+[x_4,x_6]+[x_7,x_5]   \phantom{.} \\
r_3=    [x_1,x_2]+[x_6,x_5]+[x_7,x_4]  \phantom{.} \\
r_4= [x_5,x_1]+[x_3,x_7]+[x_6,x_2]  . 
\end{array}
\;
\begin{array}{c}
r_5=  [x_1,x_4]+[x_2,x_7]+[x_3,x_6] \\
r_6=  [x_7,x_1]+[x_5,x_3]+[x_2,x_4] \\
r_7=  [x_1,x_6]+[x_4,x_3]+[x_5,x_2]  \\
\phantom{r_5=  [x_1,x_4]+[x_2,x_7]+[x_3,x_6]}
\end{array}
\end{array}
$$
Equivalently, 
$$
r_i=\sum_{m,n=1}^7 \ve^{imn} x_mx_n.
$$
Thus, if $ijk$, $ipq$, and $irs$, are distinct directed lines in the Fano plane, then 
\begin{equation}
\label{comm.reln}
r_i= [x_j,x_k] + [x_p,x_q]+[x_r,x_s].
\end{equation}
\end{prop}
\begin{pf}
Since $\mu(o_p \otimes o_q)=\sum_{i=1}^7 \ve^{ipq} o_i$,
\begin{equation}
\label{mu*}
\mu^*(x_i)= 
\sum_{p,q}  \ve^{ipq} x_p \otimes x_q.
\end{equation}
For example, 
$
o_5=o_1o_4=-o_4o_1=o_2o_7=-o_7o_2=o_3o_6=-o_6o_3, 
$
so 
$$
\mu^*(x_5)=x_1\otimes x_4 - x_4\otimes x_1 + x_2\otimes x_7 - x_7\otimes x_2 + x_3 \otimes x_6  -x_6 \otimes x_3.
$$
And so on. 
\end{pf}

 \begin{prop}
 \label{prop.relns2}
Identify $\Im {\mathbb O}_k$ with its dual via $u \leftrightarrow \langle u,-\rangle$, i.e., $x_i \equiv o_i$. 
 For each $u \in \Im {\mathbb O}_k$, there is a relation  
$$
r_u=\mu^*(\langle u,-\rangle) = \sum_{t=1}^7 \Im(o_t u) \otimes o_t 
$$
for $A$. In particular, $r_i = r_{o_i}=\mu^*(\langle o_i,-\rangle)$.
\end{prop} 
\begin{pf}
Since $\{o_1,\ldots,o_7\}$ is an orthonormal basis and $\{x_1,\ldots,x_7\}$ was defined as the dual basis the 
identification of  $\Im \mathbb O$ and  $(\Im \mathbb O)^*$ via $\langle -,-\rangle$ yields $x_i=o_i$. 

Rewriting the relation (\ref{mu*}) in terms of the $o_i$s gives 
$$
r_i  = \sum_{p,q}  \ve^{ipq} o_p \otimes o_q.
$$
However, if $i$ and $q$ are fixed, there is a unique $p$ such that $ \ve^{ipq}$ is non-zero, namely the $p$
such that $ipq$ is a line in the Fano plane, and in that case $\ve^{ipq} o_p= o_qo_i$.
Therefore
$$
r_i  = \sum_{\substack{p,q \\ p  \ne q}}  o_qo_i \otimes o_q = \sum_{q=1}^7 \Im(o_qo_i) \otimes o_q.
$$
If $u = \sum \l_i o_i$ we define $r_u =\mu^*(\sum \l_i x_i) =  \sum \l_i r_i$. 
\end{pf}

 \subsection{Alternative generators and relations for $A_\CC$}
 \label{sect.cartan}
 
In  \cite{EC1}, Cartan proved that the subgroup of $\GL(7,\CC)$ that fixes the quadratic form $t^2+u_1v_1+u_2v_2+u_3v_3$ and preserves the linear span of the seven 1-forms 
\begin{align*}
&tdu_i-u_idt+v_jdv_k- v_kdv_j \qquad (i,j,k) \; \hbox{a cyclic permutation of (1,2,3)}   \\
&tdv_i-v_idt+u_jdu_k- u_kdu_j \qquad (i,j,k) \; \hbox{a cyclic permutation of (1,2,3)}   \\
&u_1dv_1-v_1du_1+u_2dv_2-v_2du_2+u_3dv_3-v_3du_3 
\end{align*}
is a complex Lie group of type $G_2$ (see, e.g., \cite[Sect. 4]{B3}). Identify $\CC^7$ with 
the linear span of $\{dt, du_i,dv_i \; | \; i=1,2,3\}$. Let $R$ denote the subspace of alternating tensors in 
$\CC^7 \otimes \CC^7$ spanned by the exterior derivatives of the seven 1-forms. 
Then  $A_\CC$ is isomorphic to 
$T(\CC^7)/(R)$, the free algebra modulo the ideal generated by the exterior derivatives   (Proposition \ref{prop.relns/C}).   

\smallskip
 \subsubsection{Generators and relations for $A_\CC$}
 \label{sect.relns/C}
In this section we continue to assume that $\fchar k \ne 2$ and make the additional assumption 
that $k$ has a square root of $-1$ that we denote by $i$. In this case ${\mathbb O}_k$ is split. 
It is convenient to use a new basis for  ${\mathbb O}_k$,
namely
 $$
 \begin{array}{cccc}
  \begin{array}{l}
t=ix_1 \\
\phantom{t=io_1}
\end{array}
\qquad
\begin{array}{l}
u_1= x_2+ix_3 \\
v_1 =x_2-ix_3 
\end{array}
\qquad
\begin{array}{l}
u_2= x_4+ix_5 \\
v_2 = x_4-ix_5
\end{array}
\qquad
\begin{array}{l}
u_3= x_6+ix_7 \\
v_3 = x_6-ix_7.
\end{array}
\end{array}
$$
We now prove that the relations for $A_\CC$ are given by taking the exterior derivatives of the seven 1-forms listed 
in section \ref{sect.cartan}.

\begin{prop}
\label{prop.relns/C}
Let $k$ be a field of characteristic not 2 containing $i=\sqrt{-1}$. 
Then  $A_k$ is generated by $t,u_1,u_2,u_3,v_1,v_2,v_3$ subject to the relations
\begin{equation}
\label{relns/C}
\begin{array}{cc}
\begin{array}{l}
\;  [t,u_3]=[v_2,v_1]   \\
\;  [t,u_2]=[v_1,v_3]   \\
 \; [t,u_1]=[v_3,v_2]   \\
\; [u_1,v_1]+[u_2,v_2]+[u_3,v_3]=0. 
\end{array}
\qquad
\begin{array}{c}
\;  [t,v_3]=[u_1,u_2]   \\
 \; [t,v_2]=[u_3,u_1]   \\
\; [t,v_1]=[u_2,u_3]   \\
  \phantom{[u_1,v_1] .}
\end{array}
\end{array}
\end{equation}
\end{prop}
\begin{pf}
Let $r_1,\ldots,r_7$ be the relations for $A_k$ in Proposition \ref{prop.relns1}. 
Starting from the top, the four relations in the left-hand column of (\ref{relns/C}) are obtained by computing
$r_6-ir_7$, $r_5-ir_4$, $r_3-ir_2$, and $r_1$. Those in the right-hand column are obtained by computing $r_6+ir_7$, $r_5+ir_4$, and $r_3+ir_2$. 
\end{pf}

\subsubsection{Another basis for ${\mathbb O}_k$ when $\sqrt{-1} \in k$}
\label{sect.SO}
We  assume $k$ is as in section \ref{sect.relns/C}. 
Under the identification of $\Im{\mathbb O}_k$ with its dual we view the basis for 
$\Im{\mathbb O}_k^*$ in section \ref{sect.relns/C} as a basis for $\Im{\mathbb O}_k$. Thus
 $$
 \begin{array}{cccc}
 \begin{array}{l}
t=io_1 \\
\phantom{t=io_1}
\end{array}
\qquad
\begin{array}{l}
u_1= o_2+io_3 \\
v_1 =o_2-io_3 
\end{array}
\qquad
\begin{array}{l}
u_2= o_4+io_5 \\
v_2 = o_4-io_5
\end{array}
\qquad
\begin{array}{l}
u_3= o_6+io_7 \\
v_3 = o_6-io_7.
\end{array}
\end{array}
$$
The  multiplication table for ${\mathbb O}_k$ in terms of these elements is presented in Table  (\ref{mult.table.SO}). To save space  we write $w:=2(1+t) = 2 +2io_1$. 
  \begin{equation}
  \label{mult.table.SO}
\begin{array}{|r|r|r|r|r|r|r|r|}
\hline
 & t & u_1 & u_2 & u_3 & v_1 & v_2 & v_3  \\
\hline
 t & 1 & u_1 & u_2 & u_3 & -v_1 & -v_2 & -v_3  \\
u_1 & -u_1 & 0 & -2v_3 & -2v_2 & -w & 0 & 0\\
 u_2 & -u_2 & 2v_3&  0 & 2v_1 & 0 & -w & 0\\
u_3  & -u_3 & 2v_2 & -2v_1 & 0 & 0 & 0 & -w \\
v_1 & v_1 & w  & 0 & 0 & 0 & 2u_3 & -2u_2 \\
v_2 & v_2 & 0 &  w  & 0 & -2u_3 & 0  & 2u_1 \\
v_3 & v_3 & 0 & 0  & w  & 2u_2 & - 2u_1 & 0  \\
\hline
\end{array}
\end{equation}

We note that $L:=ku_1+kv_3$ is an isotropic subspace of ${\mathbb O}_k$ and $ab=0$ for all $a,b \in L$.
Its orthogonal, $L^\perp$, which is the linear span of $\{1,t,u_1,u_2,v_2,v_3\}$, 
is a subalgebra of ${\mathbb O}_k$  known as the {\sf sextonions} \cite{LM}. 

 \subsection{Constructing $A$ as a superpotential algebra}
 \label{sect.spp.alg}
 We now show that $A_\RR$ and $A_\CC$ can be constructed from 
 generic 3-forms on $\RR^7$ and $\CC^7$. Such 3-forms play a central role in the definition and study of manifolds with $G_2$ holonomy. 
 The construction of $A$ involves the now well-known procedure of forming a superpotential
 algebra, a method first invented by string theorists.

 \begin{quote}
 We will write $V=(\Im {\mathbb O}_k)^*$ for the rest section \ref{sect.spp.alg}.
  \end{quote}
  
 Let $\phi$ be the alternating trilinear form in section \ref{sect.trilin.form}.
 
 Let $\pi$ be the restriction of the   map 
$$
V^{\otimes 3} \to \wedge^3 V, \qquad  v_1 \otimes v_2 \otimes v_3 \mapsto v_1 \wedge v_2 \wedge v_3,
$$
to  the space of totally anti-symmetric tensors, i.e.,   those 3-tensors that transform 
according to the sign representation under the action of the symmetric group $S_3$,
We define the totally anti-symmetric  tensor 
\begin{equation}
\label{defn.W}
W:=\pi^{-1}(6\phi) = \sum \ve^{ijk} x_i \otimes x_j \otimes x_k.
\end{equation}
Because $\pi$ is $\GL(V)$-equivariant, $W$ is $G_2$-invariant. As Bryant points out \cite[p. 544]{B1}, the stabilizer of $W$ is slightly larger when we work over $\CC$.  

It is $W$ rather than $\phi$ that will be used to construct $A$ as a superpotential algebra.   
Consider the isomorphism
\begin{align*}
\Psi:V^{\otimes 3} \stackrel{\sim}{\longrightarrow}& \Hom( V^*, V \otimes  V), 
\\
 \Psi(v_1 \otimes v_2 \otimes v_3):= &  (\l \mapsto \l(v_1)v_2 \otimes v_3),
\end{align*}
and define
$$
\psi:=\Psi(W).
$$
Because $W$ is a non-zero $G_2$-invariant, $\psi$ is a non-zero, hence injective, homomorphism of $G_2$-modules. Therefore $ \psi(V^*)$ is the unique 7-dimensional $G_2$-submodule of $V\otimes  V$. 
Because $W$ is totally anti-symmetric, $\psi(V^*)$ consists of skew-symmetric tensors.

 \subsubsection{Superpotential algebras}
 
Let $F$ denote the free $k$-algebra on a set
of letters $L$. The {\sf cyclic partial derivative} of a word $w$ with respect to a letter $a \in L$ is
$$
{}^\circ \pd_a(w):=  \sum_{w=uav} vu.
$$ 
We extend ${}^\circ \pd_a$ to an operator on $F$ by linearity.
When  $L= \{x_1,\ldots,x_r\}$ we  usually write ${}^\circ \pd_i$ rather than ${}^\circ \pd_{x_i}$.
An element $W$ in $k\langle x_1,\ldots,x_r\rangle$ having all its homogeneous components
of degree $\ge 3$ is called a {\sf superpotential} and 
$$
{{k\langle x_1,\ldots,x_r\rangle}\over{({}^\circ \pd_i W \; | \; 1 \le i \le r)}}
$$
is called a {\sf superpotential algebra}.

\begin{prop}
\label{prop.spp.alg}
\label{prop.explicit.relns}
Let $k$ be any field and let  $W=\sum \ve^{ijk} x_i \otimes x_j \otimes x_k$. Then
  $$
 A_k \cong {{k\langle x_1,\ldots,x_7\rangle}\over{({}^\circ \pd_i W \; | \; 1 \le i \le 7)}} = {{TV}\over{(\psi(V^*))}}.
 $$
 \begin{enumerate}
 \item{}
  If  $ijk$, $ipq$, and $irs$, are directed lines in the Fano plane, then
 $$
 {}^\circ \pd_i W = \sum_{j,k} \ve^{ijk} x_j \otimes x_k =r_i = \psi(\langle x_i,-\rangle) .
 $$
  \item{}
With $r_i$ as above, 
\begin{equation}
\label{eq.W.symm}
W = \sum_{i=1}^7 x_i \otimes r_i =  \sum_{i=1}^7 r_i \otimes x_i.
\end{equation}
\end{enumerate}
\end{prop}
\begin{pf}
(1)
It follows from (\ref{defn.W}) that  ${}^\circ \pd_i W = \sum_{j,k} \ve^{ijk} x_j \otimes x_k$.
The right-hand side of this expression is the relation $r_i$ in  (\ref{mu*}). The equality 
$  \sum_{j,k} \ve^{ijk} x_j \otimes x_k =\psi(\langle x_i,-\rangle)$ follows at once from the definition of $\psi$.

(2)
This is a general fact that would appear in the proof of the non-commutative Poincar\'e Lemma stated in \cite[Prop. 1.5.13]{VG}.  Since the Poincar\'e Lemma was not 
proved in \cite{VG} we will just note that (3) follows from the computation:
\begin{align*}
W & = \sum_{i,j,k} \ve^{ijk} x_i \otimes x_j \otimes x_k \\
 & = \sum_{i=1}^7 x_i \otimes \sum_{j,k} \ve^{ijk}   x_j \otimes x_k  =  \sum_{i=1}^7 x_i \otimes r_i \\
  & = \sum_{k=1}^7 \Big(  \sum_{i,j} \ve^{ijk}   x_i \otimes x_j \Big) \otimes x_k \\
    & = \sum_{k=1}^7 \Big(  \sum_{i,j} \ve^{kij}   x_i \otimes x_j \Big) \otimes x_k  = \sum_{k=1}^7 r_k \otimes x_k
 \end{align*}
 where we have used the fact that $\ve^{ijk}=\ve^{jki}$. 
\end{pf}

\subsubsection{Remark.}
If we define  $A_k$ as a superpotential algebra, Proposition \ref{prop.A1} can be proved without having 
to know in advance that $V^{\otimes 2}$ contains a unique copy of the 7-dimensional 
irreducible representation of $G_2$ or that it is  contained in the skew-symmetric component of 
$V^{\otimes 2}$.
 
 \smallskip
 \subsubsection{Remark.}
For each $u \in V$ we define  
$$
r_u:=  \psi(\langle u,-\rangle) \in \psi(V^*).
$$
Thus $r_i=r_{x_i}$. The map $V \to \psi(V^*)$, $u \mapsto r_u$, is a $G_2$-module isomorphism 
because the maps $\psi:V^* \to \psi(V^*)$ and 
$V \to V^*$, $u \mapsto \langle u,-\rangle$, are $G_2$-module  isomorphisms.  This $r_u$ is the same
as the $r_u$ in Proposition \ref{prop.relns2}(4).

\subsection{Additional gradings on $A$}

\subsubsection{The $\ZZ_2 \times \ZZ_2 \times \ZZ_2$-grading}
An element such as $(0,1,1)$ in $\ZZ_2^3$ will be denoted by $011$, and so on. 
It is well known that ${\mathbb O}$ is a $\ZZ_2^3$-graded algebra with
\begin{align*}
&\deg 1=000, \; \phantom{{}_1} \deg o_1=001, \; \deg o_2=010, \; \deg o_3=010, \\
&\deg o_4=100, \; \deg o_5=101, \; \deg o_6=110, \; \deg o_7=111. 
\end{align*}
Since $\Im {\mathbb O}$ is a graded subspace $\Im {\mathbb O}^*$ is also a 
$\ZZ_2^3$-graded vector space. Thus $T \Im {\mathbb O}^*$ has a $\ZZ_2^3$-grading with $\deg x_i=\deg
o_i$. A simple calculation shows that $\deg r_i=\deg o_i$. It follows that $A$ becomes a $\ZZ_2^3$-graded algebra.  This is compatible with the $\ZZ$-grading on $A$ in the sense that each homogeneous component of $A$ in the $\ZZ$-grading is a $\ZZ_2^3$-graded vector space. 

We won't make use of the $\ZZ_2^3$-grading so when we speak of the degree of an element in $A$ we will mean its $\ZZ$-degree.

\subsubsection{The $\ZZ^2$-grading when ${\mathbb O}_k$ is split}

Suppose that $\sqrt{-1}$ belongs to $k$. We will use the generating set in section \ref{sect.SO} 
for $A$. There is a $\ZZ^2$-grading on $A$ given by
$$
 \begin{array}{cccc}
 \begin{array}{l}
\deg t=(0,0) \\
\phantom{t=io_1}
\end{array}
\;
\begin{array}{l}
\deg u_1= (1,0) \\
\deg v_1 = (-1,0)
\end{array}
\;
\begin{array}{l}
\deg u_2=(0,1) \\
\deg v_2 = (0,-1)
\end{array}
\;
\begin{array}{l}
\deg u_3=  (-1,-1) \\
\deg v_3 =(1,1) .
\end{array}
\end{array}
$$

\section{Automorphisms of  ${\mathbb O}_k$}
\label{sect.octs}

In this section we assume $\fchar k \ne 2$.

This section collects some results about the action of $\Aut {\mathbb O}_k$ on the sets of lines, 2-planes, and 3-planes in $\Im {\mathbb O}_k$. 

\subsection{Automorphisms of  ${\mathbb O}_k$}

Our reference for results about $\Aut {\mathbb O}_k$ is the book by Springer and Veldkamp \cite{SV}. If one is only interested in the real and complex cases, Bryant's paper \cite{B1} gives a more geometric account of the results we need.

\begin{prop}
 \cite[Cor. 1.7.3]{SV}
 \label{prop.aut.O}
Every $k$-algebra isomorphism  $\s:D \to D'$ between subalgebras of ${\mathbb O}_k$
extends to a $k$-algebra automorphism of ${\mathbb O}_k$.
\end{prop} 

\subsection{Isotropic and null subspaces of $ \Im {\mathbb O}_k$}

Let $L$ be a linear subspace of $ \Im {\mathbb O}_k$.
We write $L^2$ for the linear span of $\{uv \; | \; u,v \in L\}$. 

\subsubsection{}
Following \cite[p.149]{LM}, we say that $L$ is {\sf null}  if $L^2=0$. 
\subsubsection{}
We say $L$ is {\sf isotropic} if $\langle L,L\rangle =0$.
\subsubsection{}
A null subspace of  $ \Im {\mathbb O}_k$ is isotropic because $\langle u,v\rangle =
\Real(u\vol)$ but the converse is false when $k$ contains $\sqrt{-1}$ (Lemma \ref{lem.2-planes}).

\subsubsection{}
If $u \in  \Im {\mathbb O}_k$, then $u^2=-\langle u,u\rangle$ so a line in $ \Im {\mathbb O}_k$ is isotropic if
and only if it is null.

\begin{prop}
\label{prop.lines}
For each $\l \in k^\times/(k^\times)^2$, let $\cL_\l$ denote the set of lines $ku$ in $\Im {\mathbb O}_k$
such that $\langle u,u\rangle$ belongs to $\l$. Then $\Aut {\mathbb O}_k$ acts transitively on the set of isotropic lines and on each $\cL_\l$.
\end{prop}
\begin{pf}
Let $ku$ and $kv$ be isotropic lines in $\Im {\mathbb O}_k$. 
Then $u^2=v^2=0$ so there is an algebra isomorphism
$\s:k\oplus kv \to k \oplus ku$ such that $\s(u)=v$. By Proposition \ref{prop.aut.O}, $\s$ is the restriction of 
an automorphism of ${\mathbb O}_k$.

Suppose $ku$ and $kv$ are in $\cL_\l$. Then $\langle u,u\rangle =\a^2 \langle v,v\rangle$ for some $\a \in k^\times$. Now $k \oplus ku$ and $k \oplus kv$ are subalgebras of ${\mathbb O}_k$ and the map $\s:k \oplus ku \to k\oplus kv$ given by $\s(1)=1$ and 
$\s(u)= \a v$ is an algebra isomorphism. By Proposition \ref{prop.aut.O}, $\s$ extends to an automorphism of ${\mathbb O}_k$. 
\end{pf}

\begin{lem}
\label{lem.2-planes}
$\phantom{xxx}$
\begin{enumerate}
  \item 
  $\Aut {\mathbb O}_k$ acts transitively on the set of null 2-planes; 
  \item 
    $\Aut {\mathbb O}_k$ acts transitively on the set of non-null isotropic 2-planes.
\end{enumerate} 
\end{lem}
\begin{pf}
(1)
If $L$ and $L'$ are null 2-planes, then
$k+L$ and $k+L'$ are isomorphic subalgebras of ${\mathbb O}_k$ so there is $\s \in \Aut {\mathbb O}_k$
such that $\s(k+L)=k+L'$. It is clear that $\s(L)=L'$. This proves (2).

(2)
Let $L$ be a non-null isotropic 2-plane in $\Im {\mathbb O}_k$. Then $k+L+L^2$ is a 
subalgebra of ${\mathbb O}_k$ and is isomorphic to the exterior algebra $\L(k^2)$. If $L$ and 
$L'$ are null 2-planes in $\Im {\mathbb O}_k$ with the property that $L^2\ne 0 \ne L'^2$, then
$k+L+L^2$ and $k+L'+L'^2$ are isomorphic subalgebras of ${\mathbb O}_k$ so there is 
$\s \in \Aut {\mathbb O}_k$
such that $\s(k+L+L^2)=k+L'+L'^2$. It is clear that $\s(L)=L'$. This proves (1) and (3).
\end{pf}

\subsubsection{}
\label{ssect.root.-1}
Suppose $k$ contains a square root of $-1$, $i$ say. In the notation of section \ref{sect.SO}, $ku_1+kv_2$
is a null 2-plane and $ku_1+ku_2$ is a non-null isotropic 2-plane. In particular, there are two orbits for the action of $\Aut {\mathbb O}_k$ on the set of isotropic 2-planes.

\begin{prop}
\label{prop.2-planes}
Let $L$ be a non-isotropic  2-plane in $\Im{\mathbb O}_k$ such that $k+L$ is a subalgebra of ${\mathbb O}_k$. Then
\begin{enumerate}
  \item 
  $k+L \cong \begin{pmatrix}
  k    &  k  \\
  0    &  k
\end{pmatrix}$;
  \item 
  if $\theta:  k+L \to \begin{pmatrix}   k    &  k  \\   0    &  k \end{pmatrix}$ is an algebra isomorphism, then $\theta(L)$ is equal to the set of trace-zero matrices.
  \item 
  $\rank \langle -,-\rangle\big\vert_L =1$.
\end{enumerate}
\end{prop}
\begin{pf}
(1)
Since $\langle L,L \rangle \ne 0$, there is a basis $\{u,v\}$ for $L$ such that $\langle u,v \rangle = 1$. 

\underline{Claim:} $\Im(uv) \ne 0$. \underline{Proof:}
Suppose the claim is false. Then $uv =-1$, Since $u$ and $v$ are purely imaginary, $u^2$ and $v^2$ belong to $k.1$. Since $uv^2=-v$, $v^2 \ne 0$. 
 Hence $v$ is a unit and $v^{-1} \in kv$. 
Hence $u=(uv)v^{-1}=-v^{-1} \in kv$; this contradicts the linear independence of $u$ and $v$. The claim is therefore true. $\diamond$

Since $\Im(vu)=-\Im(uv) \ne 0$, $u$ does not commute with $v$.
Hence $k+L$ is a 3-dimensional, non-commutative, associative $k$-algebra. Up to 
isomorphism the only such algebra is the $k$-algebra of $2 \times 2$ upper triangular matrices so (1) follows.

(2)
Let $\theta$ be an isomorphism as in the statement of (2) and identify $k+L$ with its image under $\theta$.
Since $u^2$ and $v^2$ belong to $k.1$ there are $\a,\b,\c,\d \in k$ such that
$$
u=\begin{pmatrix}  \a     & \c   \\  0     & -\a \end{pmatrix} 
\qquad \hbox{and} \qquad
v=\begin{pmatrix}  \b    & \d   \\  0     & -\b \end{pmatrix} .
$$
Since $\dim(ku+kv)=2$ either $\c$ or $\d$ is non-zero. Hence
$$
0 \ne \d u-\c v = \begin{pmatrix}  \a\d-\b\c     & 0  \\  0     & \b\c-\a\d \end{pmatrix}.
$$ 
It follows that $L$ contains
$$
\begin{pmatrix}  1     & 0   \\  0     & -1 \end{pmatrix} 
\qquad \hbox{and} \qquad
\begin{pmatrix}  0   & 1  \\  0     & 0 \end{pmatrix} .
$$

(3)
We retain the notation in the proof of (2). Since $u^2=\a^2 I$ and $v^2=\b^2 I$ and $uv= 
\begin{pmatrix}  \a\b     & \a\d-\b\c   \\  0     & \a\b \end{pmatrix}$ the bilinear form $\langle -,-\rangle\big\vert_L$
is represented by the matrix 
$$
\begin{pmatrix}  -\a^2     & \a\b  \\  \a\b   & -\b^2 \end{pmatrix}.
$$
The rank of this matrix is one. 
\end{pf}

\begin{cor}
\label{cor.nd-2-planes}
$\Aut {\mathbb O}_k$ acts transitively on the set of non-isotropic 2-planes $L$ such that $k+L$ is a subalgebra of ${\mathbb O}_k$.
\end{cor}
\begin{pf}
Let $L$ and $L'$ be two such 2-planes. Let  
$$
\theta:  k+L \to \begin{pmatrix}   k    &  k  \\   0    &  k \end{pmatrix}
\qquad \hbox{and} \qquad
\theta':  k+L' \to \begin{pmatrix}   k    &  k  \\   0    &  k \end{pmatrix}
$$
 be algebra isomorphisms. 
Then $\theta(L)=\theta(L')$ so $\theta^{-1}\theta':k+L' \to k+L$ is an algebra isomorphism such that
$\theta^{-1}\theta'(L')=L$. The corollary now follows from Proposition \ref{prop.aut.O}. 
\end{pf}

In the setting of section \ref{sect.SO}, $L=kt+kv_3$ is a  non-isotropic 2-plane such that $k+L$ is a subalgebra of ${\mathbb O}_k$.

  \subsection{Quaternion subalgebras}
\label{ssect.H}
Let $\HH_k$ be the $k$-algebra with basis $1,u,v,w$ and relations $u^2=v^2=w^2=-1$,  
$uv=-vu=w$, $vw=-wv=u$, and $wu=-uw=v$. A subalgebra of ${\mathbb O}_k$ isomorphic to $\HH_k$ will be called a  {\sf quaternion} subalgebra. For example, if  $ijk$ is a directed line in the Fano plane, the linear span of $\{1, o_i,o_j,o_k\}$ is isomorphic to $\HH_k$.  

By Proposition \ref{prop.aut.O}, $\Aut {\mathbb O}_k$ acts transitively on the set of quaternion subalgebras
of ${\mathbb O}_k$.

\begin{prop}
Let $a,b \in {\mathbb O}_k$. Then $\{1,a,b\}$ is an orthonormal set if and only if $k+ka+kb+kab$ is a quaternion subalgebra of ${\mathbb O}_k$ in which $a^2=b^2=(ab)^2=-1$.
\end{prop}
\begin{pf}
($\Rightarrow$)
The hypothesis implies that $a,b \in \Im{\mathbb O}_k$ and $a^2=b^2=-1$. Also,
since $\langle a,b \rangle=0$, $ab \in \Im {\mathbb O}_k$ whence $ab=-ba$. 

E. Artin proved that any subalgebra of ${\mathbb O}_k$ generated by two elements is associative
\cite[Thm. 1.4.3]{SV}. By \cite[Prop. 1.5.1]{SV}, $k+ka+kb+kab$ is a 4-dimensional associative 
subalgebra of  ${\mathbb O}_k$. Since $ab=-ba$ it follows that $(ab)^2=-1$. It is now clear that 
 $k+ka+kb+kab$ is a quaternion subalgebra. 
 
 ($\Leftarrow$)
 Because $a^2-2\langle a,1\rangle a + \langle a,a\rangle =0$ and $a^2=-1$, $2\langle a,1\rangle a = \langle a,a\rangle -1$. But $a \notin k$, so $\langle a,1\rangle =0$ and $\langle a,a\rangle =-1$. Similarly,  
 $\langle b,1\rangle =\langle ab,1\rangle =0$ and $\langle b,b\rangle =\langle ab,ab\rangle =-1$.
 Since $b$ and $ab$ are purely imaginary, $\langle a,b\rangle = \Real(a\bol)=\Real(-ab) = - \langle ab,1\rangle =0$. Hence $\{1,a,b\}$ is an orthonormal set. 
\end{pf} 

\begin{cor}
If $\{a,b\}$ and $\{a',b'\}$ are orthonormal sets in $\Im{\mathbb O}_k$, then there is an automorphism $\s$ of ${\mathbb O}_k$ such that $\s(a)=a'$ and $\s(b)=b'$.
\end{cor}
\begin{pf}
The hypothesis implies that  $k+ka+kb+kab$  and $k+ka'+kb'+ka'b'$
are  quaternion subalgebras of ${\mathbb O}_k$ in which the squares of $a,b,ab,a',b',a'b'$ equal $-1$. 
It follows that there is an isomorphism $\s:k+ka+kb+kab\to  k+ka'+kb'+ka'b'$ such that $\s(a)=a'$ and $\s(b)=b'$. By Proposition \ref{prop.aut.O}, $\s$ extends to an automorphism of ${\mathbb O}_k$. 
\end{pf}

\subsection{Associative 3-planes and co-associative 4-planes}
\label{sect.assoc}
    A 3-plane in $\Im{\mathbb O}_k$ is {\sf associative} if it is equal to $\Im \HH$ for some quaternion subalgebra  $\HH \subset {\mathbb O}_k$. 
    Thus the associative 3-planes are the 3-planes in $\Im{\mathbb O}$ 
    having an orthonormal basis of the form $\{u,v,u v\}$. 
    It follows from Proposition \ref{prop.aut.O} that $\Aut {\mathbb O}$ acts transitively on the set of associative 3-planes.
    The associative 3-planes play a role in Proposition \ref{prop.quo.rings}.

    A 4-plane in $\Im {\mathbb O}$ is 
    {\sf co-associative} if it is orthogonal to an associative 3-plane \cite[Sect. 12.1]{J2}.  
    Because $\langle -,-\rangle$ is $\Aut {\mathbb O}$-invariant, 
    $\Aut {\mathbb O}$ acts transitively on the set of co-associative 4-planes.

\subsection{The cases $\RR$ and $\CC$}
When $k$ is $\RR$ or $\CC$ Bryant's paper \cite{B1}
gives elegant proofs of many of the foregoing results. Byant's paper also proves the facts we now state in this section.

The automorphism group of the classical octonions over $\RR$ is $G_2^c$, the compact real form of $G_2$. \

 $G_2^c$ acts irreducibly on $\Im(\mathbb O)$ and transitively on lines, 
2-planes, orthonormal pairs, 
and associative 3-planes in $\Im(\mathbb O)$ (\cite[Thm. 1, p. 539]{B1},  \cite[Sect. 2.3]{B2}, \cite[Sect. 12.1]{J2}).
As a submanifold of the Grassmannian of 3-planes in $\Im(\mathbb O)$, the set of associative 3-planes is isomorphic to $G_2/SO(4)$  \cite[Prop. 12.1.2]{J2}.

The $G_2^c$-orbit of a point $v \in \Im(\mathbb O)$ is the 6-sphere $\{x \in \Im({\mathbb O})
\; | \; \langle x,x \rangle=\langle v,v \rangle\}$ and the stabilizer of $v$ is isomorphic to $\SU(3)$, so
 $S^6 \cong G^2/SU(3)$. 
Because $SU(3)$ acts transitively on $S^5 \subset 
\RR^6$, it follows that $G_2$ acts transitively on orthonormal pairs in $\Im(\mathbb O)$  \cite[Sect. 2.3]{B2}.

More information about the octonions can be found in Baez's article \cite{Bz}.

 \section{Properties of $A$}
 
Let $k$ be a field. We will write $A_k=T(V)/(R)$ where $V=(\Im {\mathbb O}_k)^* =A_1$ and $R \subset V \otimes V$ is the linear span of the relations in Proposition \ref{prop.relns1}. 
 
In section \ref{sect.H.series} we determine a basis for $A$, use that to compute its Hilbert series,
and so deduce that $(V \otimes R) \cap (R \otimes V)$ has dimension one.
 It then follows from (\ref{eq.W.symm}) that $(V \otimes R) \cap (R \otimes V)$ is spanned by the superpotential 
 $W$. This allows us to compute $\Aut_{\hbox{\tiny{gr-alg}}} (A)$, the group of 
 graded $k$-algebra automorphisms of $A$. 
 
 In section \ref{sect.4.1} we show $A$ is an Ore extension of an algebra studied by 
 James Zhang \cite{Z}.

 \subsection{The Hilbert series for $A$}
 \label{sect.H.series}

\begin{lem}
\label{lem.basis}
Let $\cA$ be the set of words in the letters $x_1,\ldots,x_7$ that do not contain any one of 
\begin{equation}
\label{bad.words}
x_7x_6,  \quad x_7x_5, \quad   x_7x_4,  \quad x_7x_3,  \quad x_7x_2,  \quad x_7x_1,  \quad x_6x_1
\end{equation}
as a subword. Then $\cA$ is a basis for $A$.
\end{lem}
\begin{pf}
This follows from Bergman's Diamond Lemma \cite{Berg}. 
 We will use the lexicographic ordering on words with respect to $x_1<x_2< \cdots <x_7$. 
The replacement rules associated to the relations listed in  (\ref{matrix.M}) are  therefore 
\begin{align*}
x_7x_6 &= x_6x_7  + x_2x_3 - x_3x_2 + x_4x_5 - x_5x_4
\\
x_7x_5 &= x_5x_7  + x_1x_3 - x_3x_1 - x_4x_6 + x_6x_4
\\
x_7x_4 &= x_4x_7 -  x_1x_2 + x_2x_1 + x_5x_6 - x_6x_5
\\
x_7x_3 &= x_3x_7 + x_5x_1 - x_1x_5 + x_6x_2 - x_2x_6
\\
x_7x_2  &=  x_2x_7 + x_1x_4 - x_4x_1 + x_3x_6 - x_6x_3
\\
x_7x_1 & = x_1x_7 + x_4x_2 - x_2x_4 + x_3x_5 - x_5x_3
\\
x_6x_1  & =  x_1x_6 + x_5x_2 - x_2x_5 + x_4x_3 - x_3x_4.
\end{align*}
The only ambiguity is $x_7x_6x_1$. To prove the lemma we must show the ambiguity is resolvable. 
This is an easy though error-prone calculation. 
\end{pf}

  \begin{prop}
 \label{prop.Hseries}
 The Hilbert series of $A$ is $(1-7t+7t^2-t^3)^{-1}$. 
 \end{prop}
 \begin{pf}
Let $\cA$ be the basis in Lemma \ref{lem.basis}. Write $\cA_n$ for the words of length $n$ in $\cA$ and $a_n:=|\cA_n|$. 
Let $\cC$ be the words in $\cA$ that do not end in $x_6$ or $x_7$ and let $\cD$ be the  words in $\cA$ that end in $x_6$. Then
\begin{equation}
\label{eq.1}
\cA=\cA x_7 \sqcup \cC \sqcup \cD.
\end{equation}
Write $\cC_n=\cC \cap \cA_n$, $\cD_n=\cD \cap \cA_n$, $c_n=|\cC_n|$, and $d_n=|\cD_n|$. Define
$$
a(t):=\sum_{n=0}^\infty a_n t^n, \qquad  c(t):=\sum_{n=0}^\infty c_n t^n, \qquad d(t):=\sum_{n=0}^\infty d_n t^n.
$$

By Lemma \ref{lem.basis}, $\cD=(\cC \sqcup \cD)x_6$ so $d_{n+1}=c_n+d_n$ and $d(t)=tc(t)+td(t)$. Thus
$$
d(t)={{t}\over{1-t}} c(t).
$$
However,
$$
\cC =\{1\} \sqcup  \cC x_1 \sqcup \cdots  \sqcup \cC x_5 \sqcup \cD x_2  \sqcup \cdots  \sqcup \cD x_5
$$
so $c_{n+1} = 5c_n+4d_n$ for $n \ge 0$. It follows that  $c(t)-1=5tc(t) +4td(t)$ whence 
$$
(1-5t)c(t)=1+4td(t)=1+ {{4t^2}\over{1-t}} c(t).
$$
It follows that
$$
c(t)={{1-t}\over{1-6t+t^2}} \qquad \hbox{and} \qquad d(t)={{t}\over{1-6t+t^2}}.
$$
By (\ref{eq.1}), $a_{n+1}=a_n+c_{n+1}+d_{n+1}$. This implies $a(t)=(1-t)^{-1}(1-6t+t^2)^{-1}$ as claimed. 
\end{pf}

  \begin{cor}
 \label{cor.W}
Let $W=\sum_{i,j,k} \ve^{ijk} x_i \otimes x_j \otimes x_k$ be as in (\ref{defn.W}). 
Then $$(V \otimes R)\,  \cap \, (R \otimes V)=k W.$$   
 \end{cor}
 \begin{pf}
 Since 
  $$
H_A(t)= (1-7t+7t^2-t^3)^{-1} = \sum_{n=0}^\infty(7t-7t^2+t^3)^n, 
$$
 $\dim A_3=7^3-2\times 7^2 +1$. But $A_3 =V^{\otimes 3}/(V \otimes R+R \otimes V)$ so
\begin{align*}
\dim A_3   & = 7^3- \dim(V \otimes R+R \otimes V)
\\
&  = 7^3-2\dim(R \otimes V) + \dim (V \otimes R\, \cap \, R \otimes V).
\end{align*}
Therefore $\dim (V \otimes R\, \cap \, R \otimes V) =1$. But $W=\sum x_i \otimes r_i = \sum r_i \otimes x_i$
which  obviously belongs  to $(V \otimes R)\, \cap \, (R \otimes V)$.  
\end{pf}

The radius of convergence of the Hilbert series for $A$ is $3-2\sqrt{2}$ which lies between ${{1}\over{5}}$ and ${{1}\over{6}}$. We interpret this as meaning that the growth rate of $A$ lies between that of the free algebra on 5 variables and the free algebra on 6 variables.

 \subsection{Graded algebra automorphisms of $A_k$}
 
 Let $ \Aut_{\sf gr} (A)$ denote the group of $k$-algebra automorphisms of $A$ that preserve degree. 
 There is an injective homomorphism from the multiplicative group $\GG_m(k)=k^\times$ to the center of $ \Aut_{\sf gr}$   given by sending $\xi \in k^\times$ to the automorphism ``multiplication by
 $\xi^n$ on $A_n$''.

 \begin{prop}
 \label{prop.aut2}
 Suppose $k$ is $\RR$ or $\CC$. 
 There is an ``exact sequence''
 $$
 1 \to G_2 \to \Aut_{\sf gr} (A) \to \GG_m \to 1.
 $$
 \end{prop}
 \begin{pf}
 Let $R \subset V \otimes V$ denote the space of relations for $A$.
Since $A=TV/(R)$,  
$$
 \Aut_{\sf gr}(A) = \{g \in \GL(V) \; | \; g(R) \subset R\}.
$$
Since $R$ is stable under the action of $G_2$ on $V \otimes V$, the action of $G_2$ on $V =A_1$ 
extends to an action of $G_2$ as automorphisms of $A$. 
   
Now suppose that $g \in \GL(V)$ extends to a graded algebra automorphism of $A$. Since $R$ is stable under the action of $g$ so is $V \otimes R\, \cap \, R \otimes V=kW$. Define
$$
\rho:  \Aut_{\sf gr}(A)  \to \GG_m
$$
by $g\cdot W=\rho(g) W$. Since $ \Aut_{\sf gr}(A)$ contains a copy of $\GG_m$ acting by scalar multiplication on $A_1$, and  the restriction of $\rho$ to that copy of $\GG_m$  is $\xi \mapsto \xi^3$, 
$\rho$ is surjective when $k=\RR$ and when $k=\CC$. The kernel of $\rho$ is the stabilizer of $W$ in 
$\GL(V)$ which is $G_2$ (see section \ref{sect.generic}). 
This yields the sequence in the statement of the proposition. 
\end{pf}

\subsubsection{}
The subspace spanned by $x_1^2+\cdots + x_7^2$ is stable under the action of $\Aut_{\sf gr} A_k$.

\subsection{$A$ is an Ore extension of a one-relator algebra}
 \label{sect.4.1}

Let $k$ be an arbitrary field and suppose $A$ is defined by the relations in Proposition \ref{prop.relns1}.

\begin{lemma}
\label{lem.alg.B}
Let $ku$ be a non-isotropic line in $\Im{\mathbb O}_k$. Then the subalgebra of $A$ 
generated by the subspace $u^\perp$ of $\Im{\mathbb O}_k^* \equiv A_1$ is isomorphic to the algebra
\begin{equation}
 \label{eq.B.alg}
 B:={{k\langle x_1,\ldots,x_6\rangle}\over{\big([x_1,x_6]+[x_5,x_2]+[x_4,x_3]\big)}}. 
 \end{equation}
 \end{lemma}
 \begin{pf}
  Since $\Aut {\mathbb O}_k$ acts transitively on the set of non-isotropic lines it suffices to prove the lemma for
 the line $ko_7$. Since $o_7^\perp = kx_1+\cdots +kx_6$ it suffices to prove the isomorphism for the subalgebra of $A$ generated by $x_1,\ldots,x_6$. 
 
 Bergman's Diamond Lemma \cite{Berg} implies that the words in $\{x_1,\ldots,x_6\}$ 
 that do not contain $x_6x_1$ as a subword are a basis for the algebra on the right-hand side of (\ref{eq.B.alg}). 
 By Lemma \ref{lem.basis}, these words are linearly independent in $A$  so the subalgebra of $A$ 
 generated by $x_1,\ldots,x_6$ is isomorphic to $B$.
  \end{pf}

 \begin{prop}
 \label{prop.B.alg}
 Let $B$ be the  algebra in (\ref{eq.B.alg}). 
  Then 
 \begin{enumerate}
  \item 
   $\gldim B=2$;
  \item 
  $H_B(t)=(1-6t+t^2)^{-1}$;
  \item 
  $B$ is regular in the sense that 
  $  \Ext^i_B(k,B)$ is isomorphic to $k$ when $i=2$ and is zero otherwise.  
\end{enumerate}
\end{prop}
\begin{pf}
This is a very special case of the results in Zhang's paper \cite{Z}.
Specifically, (1) and (2) are given by \cite[Prop. 1.1(2)]{Z} and (3) is given by   \cite[Prop. 1.1(3)]{Z}.
\end{pf}

\begin{lem}
\label{lem.deriv}
Let $\d$ be the derivation of the free algebra $k\langle x_1,\ldots,x_6\rangle$ defined by 
  $$
 \begin{array}{ccc}
 \d(x_1)= [x_4,x_2] + [x_3,x_5]   & \quad & \d(x_4)= [x_2,x_1] + [x_5,x_6] \phantom{.} \\
  \d(x_2)= [x_1,x_4] + [x_3,x_6]   & \quad & \d(x_5)= [x_1,x_3] + [x_6,x_4] \phantom{.} \\
   \d(x_3)= [x_5,x_1] + [x_6,x_2]   & \quad & \d(x_6)= [x_2,x_3] + [x_4,x_5]. 
\end{array} 
$$
  Then   
  \begin{enumerate}
  \item 
  $\d\big([x_1,x_6]+[x_5,x_2]+[x_4,x_3]\big) =0$ and
  \item 
  $\d$ induces a derivation on the algebra $B$ defined by (\ref{eq.B.alg}).
\end{enumerate}
\end{lem}
\begin{pf}
For brevity we will write $ij$ for $[x_i,x_j]$, $[ij,k]$ for $[[x_i,x_j],x_k]$, and
 $[i,jk]$ for $[x_i,[x_j,x_k]]$. Then  
 \begin{align*}
 \d\big([x_1,x_6]+[x_5,x_2]+[x_4,x_3]\big)  = \quad & \\
 [42,6]+[35,6]+[1,23] & +[1,45] \\
+\, [13,2]+[64,2]&+[5,14]+[5,36] \\
+\, [21,3]&+[56,3]+[4,51]+[4,62]
\end{align*}
 which is zero by the Jacobi identity.
\end{pf}

We will write $\d$ for the induced derivation on $B$ and $B[X;\d]$ for the associated Ore extension.\footnote{If 
$\d$ is a derivation on an algebra 
$B$, then the Ore extension $B[X;\d]$ is, by definition,  the free left $B$-module 
with basis $1,X,X^2,\ldots$, i.e.,  $B[X;\d]=B \oplus BX \oplus BX^2 \oplus \cdots$, 
made into an associative algebra by declaring that the 
multiplication on $B[X;\d]$ is defined by $X^iX^j=X^{i+j}$ and  $Xb=bX+\d(b)$ for $b \in B$.}

\begin{prop}
\label{prop.Ore}
There is an algebra isomorphism $A \cong B [x_7;\d]$.
\end{prop}
\begin{pf}
By definition, $A$ is generated by $B$ and $x_7$
 subject to the six relations in (\ref{comm.reln}) that involve $x_7$, i.e., the relations given by $r_i=0$ for $1 \le i \le 6$. But those relations may be written as
$$
[x_7,x_i]=\d(x_i), \qquad 1 \le i \le 6
$$
so the result follows.
\end{pf}

By Proposition \ref{prop.Ore}, the Hilbert series for $A$ and $B$ are related by the formula
$$
H_A(t)=(1-t)^{-1}H_B(t)
$$
which in combination with Proposition \ref{prop.B.alg}  gives another proof that the Hilbert series of 
  of $A$ is $(1-7t+7t^2-t^3)^{-1}$.

 \subsubsection{}
 Over $\RR$, every line in $\Im {\mathbb O}$ is non-isotropic so Lemma \ref{lem.alg.B} describes every subalgebra of $A$ generated by a codimension-one subspace of $A_1$.
 
 \subsubsection{}
 If $k$ contains a square root of $-1$, say $i$, then  $\Im {\mathbb O}$ contains isotropic lines (see section \ref{ssect.root.-1}.

  By \cite[Thm. 3, p. 544]{B1}, $G_2^\CC$ acts transitively on the set of 
 null lines, and on the set of non-null lines, in $V_\CC$.

\smallskip
 \subsubsection{}
 Using Bergman's Diamond Lemma with respect to the lexicographic ordering $u_1<u_2<u_3<v_1<v_2<v_3<t=ix_1$, one easily proves that $A_\CC$ has a basis given by the words that do not contain any of the words
 $$
tu_1, \quad tu_2, \quad tu_3, \quad tv_1, \quad tv_2, \quad tv_3, \quad v_3u_3
 $$
  as a subword.
  
 \smallskip
 \subsubsection{}
 The line $\CC u_3$  is null. The subalgebra of $A_\CC$ generated by the orthogonal
 of $\CC u_3$ in $V$ is the algebra $\CC\langle t,u_1,u_2,u_3,v_1,v_2\rangle$ modulo the relations
 $$
  [t,v_2]=[u_3,u_1] , \quad [t,v_1]=[u_2,u_3], \quad  [t,u_3]=[v_2,v_1].
  $$
 This subalgebra of $A_\CC$  does  not have the good properties of  $B$.

\subsection{}
In this section $k$ is an arbitrary field.
As we noted earlier, because the defining relations for $A$ are skew-symmetric elements of 
$A_1 \otimes A_1$, $A$  is the enveloping algebra of a the Lie subalgebra $\overline{\mathfrak f}$ of $A$ generated by $A_1$.
We write $\overline{\mathfrak f}_n:=\overline{\mathfrak f} \cap A_n$.

\begin{prop}
\label{prop.Lie.alg}
The dimension of $\overline{\mathfrak f}_m$ is
$$
{{1}\over{m}} \sum_{d|m} \mu\big(\hbox{${{m}\over{d}}$}\big) \big(1+t_1^d+t_2^d\big)
$$
where $\mu$ is the M\"obius function, and $t_1$ and $t_2$ are the zeroes of $t^2-6t+1$.
\end{prop}
\begin{pf}
By the Poincar\'e-Birkhoff-Witt theorem, the Hilbert series of
$A$ is 
$$
\prod_{i=1}^\infty (1-t^i)^{-\dim \overline{\mathfrak f}_i}.
$$
But the Hilbert series of $A$ is also $(1-t)^{-1}(1-6t+t^2)^{-1}$ so
$$
\ln(1-t) + \ln(t_1-t) + \ln(t_2-t) =  \sum_{i=1}^\infty \dim \overline{\mathfrak f}_i\ln(1-t^i).
$$
Since $t_1t_2=1$, the left-hand side of the previous equation is equal to
$$
\ln(1-t)+\ln\big(1-\hbox{${{t}\over{t_1}}$}\big)   + \ln\big(1-\hbox{${{t}\over{t_2}}$}\big).
$$
It follows that
$$
\sum_{m=1}^\infty \hbox{${{1}\over{m}}$} \Big(t^m + \big(\hbox{${{t}\over{t_1}}$}\big)^m +  \big(\hbox{${{t}\over{t_2}}$}\big)^m \Big)
= \sum_{i=1}^\infty \dim \overline{\mathfrak f}_i \sum_{j=1}^\infty \hbox{${{1}\over{j}}$} t^{ij}.
$$
Equating the coefficients of $t^m$ gives
$$
 1 +   t_1^{-m}   +  t_2^{-m}
= \sum_{i|m}  i \dim \overline{\mathfrak f}_i.
$$
Because $t_1t_2=1$, $ t_1^{-m}   +  t_2^{-m} =  t_1^{m}   +  t_2^{m}$. 
The result now follows from the  M\"obius Inversion Formula.
\end{pf}

Let $v_m =  t_1^{m}   +  t_2^{m}$. Then $v_0=2$,  $v_1=6$, and  $v_m=6v_{m-1}-v_{m-2}$ for $m \ge 2$. 
The first few values of are  $v_3=34$, $198$, $1154$, $6726$, 39202, 228486, 1331714, 7761798. This is Sloane's  sequence  A003499.  It now follows from Proposition \ref{prop.Lie.alg} that the first few values of
 $\dim \overline{\mathfrak f}_m$ are
$7$, $14$, $64$, $280$, 1344, and 6496. Each $\overline{\mathfrak f}_m$ is a representation of $G_2$. The smallest dimensions of the irreducible representations of the complex Lie algebra of type $G_2$ (Sloane's sequence A104599) are 1, 7, 14, 27, 64, 77, 77, 182, 189, 273, 286, 378, 448, 714, 729, 748.

\section{Homological properties of $A$}
 \label{sect.4.2}

In this section we show $A$ is a Koszul algebra and a Calabi-Yau algebra. The methods
and notation we use are standard and we assume the reader is already familiar with these.

\begin{prop}
 \label{prop.Koszul}
Over any field, $A_k$ is a  Koszul algebra of global homological dimension 3.
\end{prop}
\begin{pf}
Let  $\xul:=(x_1,\ldots,x_7)$. The relations for $A$ can be written as a single matrix equation $M \xul^T=0$.
Here $M$ is a $7 \times 7$ matrix $M$ with entries in $A_1=V$ such that the  
$i^{\th}$ entry in $M \xul^T$, viewed as element in $TV$, is the relation $r_i$. We employ a shorthand notation for the entries in $M$: each entry im $M$ is $\pm x_j$ for some $j$ and we write $j$ for $x_j$ and $\jol$ for
$-x_j$. We also use a shorthand notation for the relations.  For example, $51.37.62$ is shorthand for the relation
$
[x_5,x_1] + [x_3,x_7] + [x_6,x_2].
$
The relations and $M$ are
  \begin{equation}
  \label{matrix.M}
 \begin{array}{ c}
r_1=   23.45.67\\
r_2=   31.46.75\\
r_3=  12.65.74 \\
r_4=  51.37.62 \\
r_5=   14.27.36\\
r_6=  71.53.24 \\
r_7=  16.43.52 
\end{array} \qquad
M = \begin{pmatrix}
 0  &\3ol  &2&\5ol&4&\7ol& 6   \\
3&0& \1ol  &\6ol&7   &4&\5ol \\
\ol2&1&0    &7&6   &\5ol&\4ol \\
5&6&\7ol  &0& \1ol  &\ol2&3 \\
\4ol&\7ol& \6ol &1&0   &3&2 \\
7&\4ol&5  &2&\3ol   &0&\1ol \\
\6ol&5&4  &\3ol&\ol2   &1&0 
\end{pmatrix} 
\end{equation}
Since $M\xul^T=0$ in $A$, $\xul M^T=0$ too. But $M$ is skew-symmetric so $\xul M=0$.
It follows from the skew-symmetry of the symbols $\ve^{ijk}$, or one can  check by hand, that
the $i^{\th}$ entry in $\xul M$ is $r_i$. 

It follows that 
\begin{equation}
\label{A.min.res}
\UseComputerModernTips
\xymatrix{
0 \ar[r] & A(-3) \ar[r]^{\cdot\xul}  & A(-2)^7 \ar[r]^{\cdot M}  &A(-1)^7 \ar[r]^<<<<<{\cdot \xul^T}  &A  \ar[r] &k \ar[r] & 0
}
\end{equation}
is a complex of left $A$-modules---entries in $A(-1)^7$ and $A(-2)^7$ are row vectors and the maps are right multiplication by the given matrices. 
Because the entries in $M\xul^T$ are a basis for the space of relations Govorov's Theorem \cite{Gov73}
shows that  the complex is exact at $A$ and $A(-1)^7$.

The map $\cdot \xul:A(-3) \to A(-2)^7$ is injective because $A$ is a domain, and a Hilbert series computation shows its image is equal to the kernel of the map 
${\cdot M}:A(-2)^7 \to A(-1)^7$. Hence (\ref{A.min.res}) is exact and therefore $A$ is Koszul.

The global homological dimension of a connected graded algebra is equal to the projective dimension
of $k$ which we have just seen is 3. 
\end{pf}

\subsection{The Koszul dual $A^!$}

The quadratic dual of $A$ is, by definition, $TA_1^*/(R^\perp)$ and because $A$ is Koszul
 is isomorphic to $\Ext^\hdot_A(k,k)$. The next result follows from the functional equation relating the Hilbert series of a Koszul algebra and its dual.

\begin{cor}
The  Hilbert series of $A^!$ is $1+7t+7t^2+t^3$. 
\end{cor}
 
 Because the symmetric algebra $S(\Im {\mathbb O}^*)$ is a quotient of $A$ by degree two elements, 
 there is a surjective homomorphism $S(\Im {\mathbb O}^*)^! \to A^!$ between their quadratic duals.  In particular, $A^!$ is a quotient of $\wedge (\Im {\mathbb O})$.

 \begin{prop}
 \label{prop.A!}
\label{prop.symm.dual}
There are identifications
$$
A_1^! = \Im {\mathbb O}, \qquad A_2^! =  \Im {\mathbb O}, \qquad A_3^! =  k,
$$
such that
\begin{enumerate}
  \item 
  the multiplication $A_1^! \times A_1^! \to A_2^!$ is the map $(u,v) \mapsto \Im(uv)$, i.e., 
  the map $\mu$ in (\ref{diag.defn.A}), and
  \item 
   the multiplication $A_1^! \times A_1^!  \times A_1^! \to A_3^!$ is the map $(u,v,w) \mapsto
   -\Real(uvw)=\psi(u,v,w)$. 
\end{enumerate}
Furthermore, 
\begin{enumerate}
  \item[(3)]
$(A^!)^* \cong A^!$ as an $A^!$-bimodule, i.e., $A^!$ is a symmetric Frobenius algebra. 
\end{enumerate}
\end{prop}
\begin{pf}
By definition, $A_1^!=\Im {\mathbb O}$.  

There are isomorphisms 
 $$
 \UseComputerModernTips
\xymatrix{
\Im {\mathbb O} \ar@{<->}[r]  & \Im {\mathbb O}^*  \ar@{<->}[r]  & R
}
$$
such that 
$$
 \UseComputerModernTips
\xymatrix{
w \ar@{<->}[r]  & \langle w,-\rangle  \ar@{<->}[r]  & r_w = \mu^*\big(\langle w,-\rangle \big) 
= \mu^*\big(\langle -,w\rangle \big) 
}
$$
and
$$
A_2^! = \frac{A_1^! \otimes A_1^!}{R^\perp} =  \frac{\Im {\mathbb O}\otimes \Im {\mathbb O}}{R^\perp} \cong R^*.
$$
so $A_2^!$ naturally identifies with $\Im {\mathbb O}$.  

With this identification, the multiplication $A_1^! \times A_1^! \to A_2^!$ corresponds to a bilinear 
map $\Im {\mathbb O}_k  \times \Im {\mathbb O}_k  \to \Im {\mathbb O}_k$.   To avoid confusion with multiplication in ${\mathbb O}_k$ we write $*$ for the product in $A^!$. 

(1)
Let $u,v \in A_1^!=\Im {\mathbb O}_k$. The value of $u*v \in A_2^!$  on $r_w \in R$ is   
\begin{align*}
(u*v)(r_w) & =(u \otimes v)(r_w)
\\
& =(u \otimes v)\Big( \mu^*\big(\langle -,w \rangle \big)\Big) 
\\
& = \langle \mu(u \otimes v),w\rangle
\\
&= \langle \Im(uv),w\rangle
\end{align*}
so $u*v=\Im(uv)=\mu(u \otimes v)$. 
 
(2)
It is clear that $A_3^!$ 
naturally identifies with $(R\otimes A_1  \cap A_1 \otimes R)^*$. In our case, $R\otimes A_1  \cap A_1 \otimes R$ is spanned by $\sum_{i=1}^7 x_i \otimes r_i=\sum_{i=1}^7 r_i \otimes x_i$. 
If $u,v,w \in \Im {\mathbb O}$, then
\begin{align*}
\big(u*v*w\big)\Big(\sum_{i=1}^7 x_i \otimes r_i\Big) & = \sum_{i=1}^7x_i(u). (v*w)(r_i)
\\
&= \sum_{i=1}^7 x_i(u) \langle o_i,\Im(vw) \rangle 
\\
&=    \langle u,\Im(vw) \rangle 
\\
&=-\Real(uvw)
\\
& = \psi(u,v,w)
\end{align*}
This completes the proof.

(3)
The condition that $A^!$ be a symmetric Frobenius algebra is equivalent to the existence of a non-degenerate 
bilinear form $f:A^! \times A^! \to k$ such that $f(u*v,w)=f(u,v*w)$ and $f(u,v)=f(v,u)$ for all $u,v,w \in A^!$. 
We identify $k$ with $A_3^!$ as in the statement of the proposition and define $f$ as follows:
 $f(A^!_m,A_n^!)=0$ if $m+n \ne 3$ and $f(a,b)=a*b$ if $a \in A_m^!$ and $b \in A_{3-m}^!$. Because elements of even degree belong to the center of $A^!$, $f$ is symmetric. It is clear from the definition that 
 $f(a*b,c)=f(a,b*c)$.

 Let $u \in A_1^!-\{0\}$. 
Then $\Im(uo_i) \ne 0$ for some $i$, and $o_i=o_jo_k$ for some  $j$ and $k$ so
 $f(u,o_j*o_k) \ne 0$. Hence $f$ is non-degenerate.
\end{pf}

Part (3) of Proposition \ref{prop.A!} is a special case of the next result.

\begin{prop}
Let $A$ be a connected graded algebra that is generated in degree one. If $\gldim A=3$  and $A$ 
is  Koszul and Gorenstein, then $A^!$ is a symmetric Frobenius algebra.
\end{prop}
\begin{pf}
The Gorenstein and Koszul hypotheses imply that $\dim A_3^! =1$. It therefore suffices toshow there is a  non-degenerate bilinear form $f:A^! \times A^! \to A_3^!$ such that 
$f(ab,c)=f(a,bc)$ and $f(a,b)=f(b,a)$ for all $a,b,c \in A^!$. 

Let  $a \in A^!_m$ and $b \in A_n^!$. We define $f(a,b)=ab$ if $m+n=3$  and $f(a,b)=0$ otherwise. 
It is clear that $f(ab,c)=f(a,bc)$.

Let $R \subset A_1 \otimes A_1$ be the space of relations for $A$. Then $A_3^!$ naturally
identifies with the dual of $A_1 \otimes R \cap R \otimes A_1$. 
 Because $A$ is a superpotential algebra, the Poincar\'e lemma tells us the following: if $x_i$ is a basis
 for $A_1$ there is a corresponding basis $r_i$ for $R$ such that $\sum x_i \otimes r_i=\sum r_i \otimes x_i$. 
 Furthermore, the Gorenstein hypothesis implies that $\dim A_1=\dim R$ so all $r_i$ are non-zero.
 Hence, if $a \in A_1^!$ and $b \in A_2^!$, then 
 $$
 ab\Big(\sum x_i \otimes r_i\Big)=\sum a(x_i) b(r_i) = \sum b(r_i)a(x_i) = ba \Big(\sum r_i \otimes x_i\Big);
$$
thus $ba=ab$ and $f(a,b)=f(b,a)$. If $a \in A_0^!$ and $b \in A_3^!$, then $ab=ba$. Hence $f$ is symmetric.

Let $a \in A_1^!-\{0\}$. Then we can choose a basis $\{x_i\}_i$ for $A_1$ such that $a(x_1) \ne 0$ 
and $a(x_i)=0$ for all other $i$. Since $r_1 \ne 0$ there is $b \in A_2^!$ such that $b(r_1) \ne 0$;
it follows that 
$$
f(a,b)\Big(\sum x_i \otimes r_i\Big)=\sum a(x_i) b(r_i) =a(x_1)b(r_1) \ne 0
$$
thus showing that $f$ is non-degenerate.
\end{pf}

When the base field is $\RR$ or $\CC$, let $\wedge^2_{14}V^*$ denote the 14-dimensional irreducible
$G_2$-submodule of $\wedge^2 V^*$; then 
$$
A^! \cong {{\wedge V^*}\over{(\wedge^2_{14}V^*)}}.
$$
Furthermore, $A^!_3$ is, in effect, the span of ``the unique'' $G_2$-invariant 3-form in $\wedge^3 V^*$.
The multiplication in $A^!$ gives non-degenerate $G_2$-equivariant pairings $A_1^! \times A_2^! \to A_3^!$
and  $A_2^! \times A_1^! \to A_3^!$.

 Let $\{\xi^i \; | \;1 \le i \le 7\}$ be the basis for $A_1^!$ dual to $\{x_i \; | \; 1 \le i \le 7\}$
and define 
$$
{\bf e}:= \sum_{i=1}^7 x_i \otimes \xi^i  \qquad \hbox{and} \qquad \overline{\bf e}:= \sum_{i=1}^7 \xi^i \otimes x_i.
$$
Then the resolution (\ref{A.min.res}) is isomorphic to the Koszul complex $K_\ldot(A)$ which can be written more formally as  
\begin{equation}
\label{A.min.res}
\UseComputerModernTips
\xymatrix{
0 \ar[r] & A \otimes W \ar[r]^{d_3}  & A \otimes R   \ar[r]^{d_2}  &A\otimes A_1 \ar[r]^{d_1} &A  \ar[r] &k \ar[r] & 0
}
\end{equation}
where 
$$
d_n(a \otimes u):= {\bf e}\cdot (a \otimes u) = \sum_{i \in I} a x_i \otimes \xi^i \cdot u.
$$

\begin{cor}
\label{cor.reg}
The trivial $A$-module ${}_Ak$ has the property 
\begin{equation}
\label{eq.reg.alg}
\Ext_A^i(k,A) \cong \begin{cases}
k(3) & \text{ if $i=3$, and}
\\
0 &  \text{ if $i\ne 3$.}
\end{cases}
\end{equation}
\end{cor}
\begin{pf}
Because $K_\ldot(A)$ is a projective resolution of ${}_Ak$, $\Ext^\ldot(k,A)$ is the homology of  the complex $\Hom_A(K_\ldot(A),A)$. However,  there are isomorphisms of complexes of right $A$-modules
$$
\Hom_A(K_\ldot(A),A) \cong  (A_{\ldot}^! \otimes A, \cdot \overline{\bf e}) \cong ((A^! )^*_{3-\ldot}\otimes A, \cdot \overline{\bf e})
$$
where the second isomorphism follows from the fact that $A^!$ is a symmetric algebra.
Thus $\Hom_A(K_\ldot(A),A)$ is isomorphic to the Koszul complex for $A$ on the other side 
which, by essentially the same argument as in 
Proposition \ref{prop.Koszul}, is a projective resolution for $k_A$.  This proves (\ref{eq.reg.alg}).
\end{pf}

\begin{cor}
\label{cor.not.noeth}
The algebra $A$ is neither left nor right noetherian.
\end{cor}
\begin{pf}
Let $I$ be the ideal $(x_5,x_6,x_7)$. Then $A/I$ is isomorphic to the polynomial ring  $k[x_1, x_2,x_3,x_4]$. 
 If $A$ were left or right noetherian, then $A/I$
would have a finite free resolution in which each term is a finitely generated free $A$-module. It would follow that the Hilbert series of $A/I$ is $p(t)H_A(t)$ for some $p(t) \in \ZZ[t]$. But that is impossible because the Hilbert series of $A/I$ is $(1-t)^{-4}$. 
\end{pf}

\begin{thm}
\label{thm.CY}
$A$ is a  Calabi-Yau algebra of dimension 3.
\end{thm}
\begin{pf}
Bocklandt's criterion \cite[Thm. 4.2]{Bo} tells us it suffices to produce an $A$-bimodule resolution
\begin{equation}
\label{CY.res}
\UseComputerModernTips
\xymatrix{
0 \ar[r] & P_3 \ar[r]^{d_3}  & P_2  \ar[r]^{d_2}  & P_1  \ar[r]^{d_1} & P_0  \ar[r] &A \ar[r] & 0,
}
\end{equation}
having the property that the $P_i$s are finitely generated projective bimodules and
there are bimodule isomorphisms $\a_i$ fitting into a commutative diagram
$$
\UseComputerModernTips
\xymatrix{
  P_3 \ar[rr]^{d_3} \ar[d]^{\a_3}   && P_2  \ar[rr]^{d_2} \ar[d]^{\a_2}   && P_1  \ar[rr]^{d_1} \ar[d]^{\a_1}   && P_0  \ar[d]^{\a_0}  
\\
 P_0^\vee \ar[rr]_{-d_1^\vee}  && P_1^\vee  \ar[rr]_{-d_2^\vee}  && P_2^\vee  \ar[rr]_{-d_3^\vee} && P_3^\vee 
}
$$
where $P^\vee =\Hom_{\sf Bimod}(P,A \otimes A)$. In the above $\Hom_{\sf Bimod}$ 
is taken with respect to the {\it outer} bimodule structure on $A \otimes A$, namely $x(a \otimes b)y=xa \otimes by$, and $P^\vee$ is viewed
as a bimodule through the surviving {\it inner} bimodule structure on $A \otimes A$, namely $x*(a \otimes b)
*y=ay \otimes xb$.

Because $A$ is Koszul, Proposition 3.1 of Van den Bergh's paper \cite{VdB0} provides a 
projective resolution of $A$ as an $A$-bimodule of the form (\ref{CY.res}) in which
$$
P_n :=A \otimes (A_n^!)^* \otimes A 
$$
and
\begin{align*}
d_n(a \otimes t \otimes  b):=& \sum ax_i \otimes t \xi^i \otimes b +(-1)^n a \otimes \xi^i t \otimes x_i b
\\
= & (a \otimes  t \otimes  b)({\bf e}\otimes 1) +(-1)^n (1 \otimes \overline{\bf e})(a \otimes  t \otimes  b).
\end{align*}
Now
$$
P_n^\vee := \Hom_{A-\Bimod}(P_n,A\otimes A) = A \otimes A_n^! \otimes A
$$
and 
\begin{align*}
d_n^\vee(u \otimes \tau \otimes v) =&
\sum u \otimes  \xi^i \tau \otimes x_i v + (-1)^n ux_i \otimes \tau \xi^i  \otimes v
\\
= & (1 \otimes \overline{\bf e})(u \otimes \tau \otimes v) +  (-1)^n (u \otimes \tau \otimes v)({\bf e} \otimes 1).
\end{align*}

 By Proposition \ref{prop.symm.dual}, there is an isomorphism $\beta:(A^!)^* \to A^!$ of $A^!$-bimodules
 with components $\beta_n: (A^!_n)^* \to A^!_{3-n}$. Let 
 $$
 \a_n=
 \begin{cases}
	 -( \id_A \otimes \b_n \otimes \id_A) & \text{if $n \equiv 0,1\,$(mod 4)}
	 \\
	 + ( \id_A \otimes \b_n \otimes \id_A) & \text{if $n \equiv 2,3\,$(mod 4).}
\end{cases}
$$
Thus $\a_n:  P_n = A \otimes (A^!_n)^*\otimes A  \to P_{3-n}^\vee = A \otimes A^!_{3-n} \otimes A$
is an $A$-bimodule isomorphism. An easy calculation shows that the rectangle
$$
\UseComputerModernTips
\xymatrix{
  P_n \ar[rr]^{d_n} \ar[d]_{\a_n}   && P_{n-1}  \ar[d]^{\a_{n-1}}    
  \\
 P_{3-n}^\vee \ar[rr]_{-d_{4-n}^\vee}  && P_{4-n}^\vee  
 }
$$
 commutes, thus showing that Bocklandt's criterion holds.
\end{pf}

The following ``Poincar\'e duality'' between Hochschild homology and cohomology of $A$-bimodules
follows  from \cite[Thm. 1]{VdB2}.

\begin{cor}
If $M$ is an $A$-bimodule, them $HH_\ldot(M)  \cong HH^{3-\udot}(M)$.
\end{cor}

\begin{cor}
If $N$ is a left $A$-module, then $\Ext_A^\udot(k,N) \cong \Tor_{3-\ldot}^A(k,N)$.
\end{cor}
\begin{pf}
Since $A$ is a Hopf algebra, Theorem \ref{thm.CY} allows us to apply Lemma 2.1 in \cite{Far} which gives the result. 
\end{pf}

\section{Some representations of $A_\RR$ and $A_\CC$}
 
At present we don't know much about finite-dimensional representations of $A$. 
We make two general remarks. 

First,  since $A_k$ is an  enveloping algebra of a Lie algebra, if $M$ and $N$ are $A$-modules so is $M \otimes_k N$. When $k$ is algebraically closed there are homomorphisms
from $A$ onto the enveloping algebra of $\fsl(2,k)$ (see below) so one may use diagonal
maps $A \to U(\fsl(2))^{\otimes n}$ to obtain $A$-modules.

Second, the next result shows that the representation theory of $A$ is related to a linear algebra problem 
 over ${\mathbb O}$. 

\begin{prop}
Let $X_1,\ldots,X_7 \in M_n(k)$ and define
$$
X:=\sum_{i=1}^7 X_i o_i \in  M_n(\Im {\mathbb O}_k) \subset M_n({\mathbb O}_k).
$$
The assignment $x_i \mapsto X_i$ gives $k^n$ an $A$-module structure if and only if $X^2\in M_n(k)$,
i.e., $\Im(X^2)=0$. 
\end{prop}
\begin{pf}
Since 
$$ 
X^2 =  \sum_{\ell=1}^7 \Big(  \sum_{i,j=1}^7  \ve^{\ell ij} X_iX_j \Big)o_\ell - \sum_{\ell=1}^7 X_\ell^2
$$ 
the result follows from the fact that the defining relations for $A$ are $r_\ell=\sum_{i,j=1}^7  \ve^{\ell ij} x_ix_j$ for $\ell =1,\ldots,7$. 
\end{pf}

Thus, $n$-dimensional $A$-modules annihilated by  $x_1^2+\cdots+x_7^2$ correspond to elements $X \in M_n({\mathbb O})$ such that $X^2=0$. If $M$ and $M'$ are the $A$-modules corresponding to 
matrices $X,X' \in M_n({\mathbb O})$, then $M \cong M'$ if and only if $X'=gXg^{-1}$ for some $g \in \GL_n(k)$, so the representation theory of $A$ is not by any means equivalent to a linear algebra problem over
 ${\mathbb O}$ in that only very special kinds of changes of basis are allowed. 
 
 \subsection{Free algebras as quotients of $A_k$}
 
\begin{prop}
\label{prop.free.quo}
Let $F_k$ be the free $k$-algebra on two degree-one generators. 
\begin{enumerate}
  \item 
 If $\sqrt{-1} \in k$, then $F_k$ is a quotient of $A_k$ by a graded ideal. 
  \item 
 $F_\RR$ is not a quotient of $A_\RR$ by a graded ideal. 
\end{enumerate}
\end{prop}
\begin{pf}
(1)
Using the basis and relations in section \ref{sect.relns/C}, one sees that $A/(t,u_1,v_1,u_2,v_3)$ 
is isomorphic to $F_k$.  

(2)
Let $\varphi:A_\RR \to F_\RR$ be a homomorphism of graded $\RR$-algebras. Since $\dim_k F_1=2$, $\varphi$ vanishes on a 5-plane in $V$. Since $\Aut A_\RR$ acts transitively on 5-planes in $A_1$  
we can assume that 
$\ker \varphi$ contains $x_1,\ldots, x_5$. But $A/(x_1,\ldots,x_5)$ is the polynomial 
ring on two variables so $\varphi$ is not surjective. In fact, the image of $\varphi$ is isomorphic to either $\RR$ or the polynomial ring $\RR[X]$. 
\end{pf}

\subsubsection{Remark}
The image of $A_\RR$ in 
the free $\CC$-algebra $A_\CC/(t,u_1,v_1,u_2,v_3)$ 
is isomorphic to $\RR \langle x_4,x_5,x_6,x_7 \rangle$ modulo the relations
$$
\begin{array}{cc}
\begin{array}{l}
x_7x_4=-x_6x_5 \\
x_4x_7=-x_5x_6 \\
x_7x_5= \phantom{-}x_6x_4 \\
x_5x_7=\phantom{-} x_4x_6 
\end{array}
\qquad
\begin{array}{l}
x_4^2+x_5^2 =0  \\
x_6^2+x_7^2=0  \\
x_4x_5-x_5x_4=x_7x_6-x_6x_7. \\  \phantom{x} 
\end{array}
\end{array}
$$

\begin{cor}
For every integer $n \ge 1$, $A_\CC$ has an irreducible representation of dimension $n$, and 
there is a surjective $\CC$-algebra homomorphism $A_\CC \to M_{n}(\CC)$.
\end{cor}
\begin{pf}
The Lie algebra $\fsl(2,\CC)$   is generated by two elements. Its 
enveloping algebra  $U(\fsl(2,\CC))$ is therefore generated by two elements.
Hence $U(\fsl(2,\CC))$ is a quotient of $A_\CC$ by Proposition \ref{prop.free.quo}(1). Since 
 $\fsl(2,\CC)$ has an irreducible representation of dimension $n$ for every integer $n \ge 1$
 so does $A_\CC$. 
\end{pf}

\begin{cor}
\label{cor.MC.quo}
For every integer $n \ge 0$, there is a surjective $\RR$-algebra homomorphism $\varphi:A_\RR 
\to M_{n+1}(\CC)$.
\end{cor}
\begin{pf}
Let $\cD$ denote the ring of differential operators on the polynomial ring $\CC[x,y]$.
There is an $\RR$-algebra homomorphism $A_\RR \to \cD$ given by 
$$
x_1 \mapsto 0, \;x_2 \mapsto 0, \;x_3 \mapsto 0, \;x_4 \mapsto -ix\pd_y, \;x_5 \mapsto x\pd_y, \;x_6 \mapsto iy\pd_x, \;x_7 \mapsto y\pd_x.
$$

We will now show that each homogeneous component  $\CC[x,y]_n$ is a simple $A_\RR$-module. 
Let $0 \ne f \in \CC[x,y]_n$. Let $M$ be the $A_\RR$-submodule of $\CC[x,y]_n$ generated by $f$.
A suitable power of $x_5$ acts on $f$ to  produce a non-zero element $\l x^n$ in $M$, $\l \in \CC$.  
Hence $M$ contains $\RR x_4x_7\cdot \l x^n + \RR \l x^n=\CC x^n$. 
It follows that $M$ contains $ x_7^p\cdot \CC x^n=\CC x^{n-p}y^p$ for all $0 \le p \le n$, i.e., $M=\CC[x,y]_n$,
as claimed. 

The endomorphism ring of $\CC[x,y]_n$ as an $A_\RR$-module is $\CC$, so the result follows. 
\end{pf}

\subsection{The Lie algebra $\fgol$ and the Connes-Dubois-Violette algebra}
\label{ssect.A+}
In this section $k$ is an arbitrary field and we write $A$ for $A_k$.

As before, $\overline{\mathfrak f}$ is the Lie subalgebra of $A$ generated by $A_1$. We write
 $\overline{\mathfrak f}_n:= \overline{\mathfrak f} \cap A_n$.  
Let $\fgol$ be the free Lie algebra over $k$ on $t_0,t_1,t_2,t_3$ modulo the relations 
$$
[t_0,t_i]+[t_j,t_k]=0
$$
where $(i,j,k)$ runs over the cyclic permutations of $(1,2,3)$. We  view $\fgol$ as an $\NN$-graded Lie algebra with $\fgol_1= k t_0+kt_1+kt_2+k t_3$.

Connes and Dubois-Violette write $\cA^{(+)}$ for the enveloping algebra of $\fgol$  
and call it the {\sf quadratic self-duality algebra}   \cite[Sect. 4]{CDV}.  
Clearly, 
$$
{{A}\over{(x_1,x_2,x_3)}} \cong \cA^{(+)} = U(\fgol).
$$
The subalgebra of $U(\fgol)$ generated by any three of the $t_i$s is a free algebra,
and that $U(\fgol)$ is an Ore extension of that subalgebra. Furthermore, by  \cite{CDV}, 
$U(\fgol)$ 
 \begin{itemize}
  \item 
is a Koszul algebra of global dimension two;
  \item 
has Hochschild dimension two;
  \item 
has Hilbert series $(1-3t)^{-1}(1-t)^{-1}$, and is therefore not Gorenstein.
\end{itemize}
The following observation is easily proved. 

\begin{lem}
\label{lem.CDV}
Let $z$ be a non-zero element in $U(\fgol)_1$ and let $J$ be the ideal in $U(\fgol)$ generated by the elements
$[z,t_i]$, $4 \le i \le 7$. Then $U(\fgol)/J$ is a polynomial ring in four variables. 
\end{lem}

\subsection{Representations of $A_\RR$}

\begin{lemma}
\label{lem.Lie.hom}
Let $\overline{\mathfrak f}=\overline{\mathfrak f}_\RR$ be as above. Then
\begin{enumerate}
  \item 
  if $\fa$ is the ideal in $\overline{\mathfrak f}$ generated by an associative 3-plane in $\overline{\mathfrak f}_1$, then $\overline{\mathfrak f}/\fa \cong \fgol$;
  \item 
  if $\fb$ is an ideal in $\fgol$ such that $\fb_1 \ne 0$, then $\fgol/\fb$ is abelian;
  \item 
  if $\fc$ is an ideal in $\overline{\mathfrak f}$ such that $\dim \fc_1 \ge 4$, then $\overline{\mathfrak f}/\fc$ is abelian;
  \item{}
  Neither $\fsl(2,\RR)$, nor $\fso(3,\RR)$, nor the 3-dimensional Heisenberg Lie algebra, nor the 2-dimensional   non-abelian Lie algebra, is a quotient of $\overline{\mathfrak f}$;
  \item{}
  $\fso(3,1)$ is a quotient of $\fgol$ and hence of $\overline{\mathfrak f}$.
\end{enumerate}
\end{lemma}
\begin{pf}
(1)
Since $\Aut A_\RR$ acts transitively on the associative 3-planes we can and will assume that 
$\fa=\langle x_1,x_2,x_3\rangle$. The relations $r_4$, $r_5$, $r_6$, and $r_7$, become vacuous in $\overline{\mathfrak f}/\fa$,
and the relations $r_1$, $r_2$, and $r_3$, become those for $\fgol$ under the homomorphism $x_i \mapsto t_{i-4}$ for $i=4,5,6,7$. 

(2)
By assumption, $\fb$ contains an element 
$$
t:=  a t_0+bt_1+ct_2 +dt_3
$$
for some  $(a,b,c,d) \in \RR^4-\{0\}$. Taking commutators of $t$ 
with $t_0$, $t_1$, $t_2$, and $t_3$, it follows that 
$\fb$ contains the elements
\begin{align*}
b [t_0,t_1]+c[t_0,t_2]+d[t_0,t_3]
\\
 a [t_1,t_0]+c[t_1,t_2] + d[t_1,t_3] &= -a[t_0,t_1] -c[t_0,t_3]+d[t_0,t_2]
\\
  a [t_2,t_0]+b[t_2,t_1]+  d[t_2,t_3]  &= -a[t_0,t_2] +b[t_0,t_3]-d[t_0,t_1]
\\
a [t_3,t_0]+b[t_3,t_1]+c[t_3,t_2]   &= -a[t_0,t_3] -b [t_0,t_2]+c[t_0,t_1]. 
\end{align*}
Equivalently, the image of $t_0$ in $\fgol/\fb$  commutes with the entries in the matrix
$$
\begin{pmatrix}
 b     &   c & d \\
  a    &   -d & c   \\
  d & a & -b \\
c & -b & -a  
\end{pmatrix}
\begin{pmatrix}
t_1 \\ t_2 \\ t_3
\end{pmatrix}.
$$
Over $\RR$, the rank of this $4 \times 3$ matrix is $3$ so the image of $t_0$ in $\fgol/\fb$ commutes with
the images of $t_1$, $t_2$, and $t_3$. It now follows from the relations $[t_0,t_i]+[t_j,t_k]=0$ that 
$\fgol/\fb$ is abelian.   

(3)
Since $\Aut A_\RR$ acts transitively on the set of 2-planes  in $\Im{\mathbb O}$   
we can and will assume that $\fc$ contains $\langle x_1,x_2\rangle$. Inspecting the relations, one sees that $\overline{\mathfrak f}/\langle x_1,x_2\rangle \cong \fgol \oplus \RR$ where $\RR$ is the center of  $\fgol \oplus \RR$. 
Let $\fcol$ denote the image of $\fc$ in $\overline{\mathfrak f}/\langle x_1,x_2\rangle$.
Then $\fcol \cap \fgol \ne 0$, so it follows from (2) that $(\fgol \oplus \RR)/\fcol$ is abelian.  

(4)
A Lie algebra homomorphism from $\overline{\mathfrak f}$ to any of the indicated Lie algebras must vanish on a 4-dimensional
subspace of $\overline{\mathfrak f}_1$ so the result follows from (3).

(5)
The Lie algebra $\fso(3,1)$ has a basis consisting of the $4 \times 4 $ matrices 
 \begin{equation}
 \label{eq.so31.basis}
A_i:= e_{jk}-e_{kj}  \qquad \hbox{and} \qquad B_i:=e_{i4}+e_{4i} 
 \end{equation}
where $(i,j,k)$ runs over the cyclic permutations of $(1,2,3)$.
If $(i,j,k)$ is a cyclic permutation of $(1,2,3)$,  then
$$
[A_i,A_j]=-A_k, \; [A_i,B_i]=0, \; [A_i,B_j]=  [B_i,A_j]=-B_k,  \; [B_i,B_j]=A_k.
$$
The maps $x_1,x_2,x_3 \mapsto 0$ and 
$$
x_4 \mapsto t_0 \mapsto A_1, \quad  x_5 \mapsto t_1 \mapsto A_2, \quad
 x_6 \mapsto   t_2 \mapsto B_1, \quad  x_7 \mapsto  t_3 \mapsto B_2.
$$
give surjective Lie algebra homomorphisms $\overline{\mathfrak f} \twoheadrightarrow \fgol \twoheadrightarrow \fso(3,1)$.
\end{pf}

 \begin{prop}
$\phantom{xxx}$
\begin{enumerate}
  \item{}
  The  image of an $\RR$-algebra homomorphism from $A_\RR$ to either $M_2(\RR)$, or $\HH$, or $\begin{pmatrix}
\RR      &    \RR \\
      0 & \RR 
\end{pmatrix}$, is isomorphic to either $\RR$, $\CC$, or $\RR[\ve]/(\ve^2)$.
   \item 
  Neither $M_2(\RR)$, nor $\HH$, nor $\begin{pmatrix}
\RR      &    \RR \\
      0 & \RR 
\end{pmatrix}$, is a quotient of $A_\RR$
  \item 
  $M_4(\RR)$ is a quotient of $A_\RR$.
  \end{enumerate}
\end{prop}
\begin{pf}
(1)
Let $\phi$ be an $\RR$-algebra homomorphism from $A_\RR$ to $M_2(\RR)$ or $\HH$.
Then $\phi$ restricts to a Lie algebra homomorphism from $\overline{\mathfrak f}$ to $\RR \oplus \fsl(2,\RR)$ or
$\RR \oplus \fso(3,\RR)$. If $\dim \phi(\overline{\mathfrak f})=4$, then $\fsl(2,\RR)$ or
$\fso(3,\RR)$ would be a quotient of $\overline{\mathfrak f}$ thereby contradicting Lemma \ref{lem.Lie.hom}(4).
Hence $\dim\phi(\overline{\mathfrak f}) \le 3$, Now Lemma \ref{lem.Lie.hom}(3) implies that $\phi(\overline{\mathfrak f})$ is an abelian Lie 
algebra and hence that $\phi(A_\RR)$ is a commutative ring. This also proves (2)

(3)
The proof of Corollary \ref{cor.MC.quo} shows that $\fsl(2,\CC)$ is a quotient of the real Lie algebra $\overline{\mathfrak f}$. 
As a real Lie algebra $\fsl(2,\CC)$ is isomorphic to $\fso(3,1)$ so $\fso(3,1)$ is a quotient of $\overline{\mathfrak f}$ and therefore the enveloping algebra $U(\fso(3,1))$ is a quotient of $A_\RR$. 
It is easy to see that the associative subalgebra of $M_4(\RR)$ generated by the basis for 
$\fso(3,1)$ listed in (\ref{eq.so31.basis}) is $M_4(\RR)$. Hence $M_4(\RR)$ is a quotient of $U(\fso(3,1))$ 
and therefore a quotient of $A_\RR$. 
\end{pf}

\begin{cor}
\label{cor.no.exts}
If $M$ and $N$ are non-isomorphic 1-dimensional representations of $A_\RR$, then $\Ext^1_A(M,N)=0$.
\end{cor}

\subsubsection{Remark}
I do not know if $M_3(\RR)$ is a quotient of $A_\RR$.

\subsubsection{Remark}
Because the complexification of $\fso(3,1)$ is isomorphic to $\fsl(2,\CC) \times \fsl(2,\CC)$, there is a surjective homomorphism
$$
\overline{\mathfrak f}_\CC \twoheadrightarrow  \fsl(2,\CC) \times \fsl(2,\CC).
$$

\subsubsection{Remark}
Since the polynomial ring $\RR[t_1,\ldots,t_4]$ is a quotient of $A_\RR$, so is the polynomial ring 
$\CC[t_1,t_2,t_3]$.

\section{Ideals and left ideals of $A_\RR$ generated by subspaces of $A_1$}
 
Since the relations for $A_k$ are skew symmetric there is a surjective $k$-algebra homomorphism $A_k \to S(\Im {\mathbb O}_k^*)$, the polynomial ring on seven variables. 

Let $\cA^{(+)}$ be the enveloping algebra $U(\fgol)$ in section \ref{ssect.A+}. 
 
 \begin{prop}
 \label{prop.quo.rings}
Suppose the base field is $\RR$. 
 \begin{enumerate}
  \item 
  If $L$ is a 2-plane in $A_1$, then  $A/(L) \cong \cA^{(+)} \otimes \RR[t]$ where $t$ is a central indeterminate. 
  \item{}
 If $P$ is a 3-plane, then  
  $$
  {{A}\over{(P) }} \cong 
    \begin{cases}
   & \text{$ \cA^{(+)}$ if $P$ is associative (section \ref{sect.assoc}), and}
  \\
  & \text{a polynomial ring on four variables otherwise.}
  \end{cases}
  $$
  \item{}
  if $H$ is a $d$-plane in $V$ and $d \ge 4$, then $A/(H)$ is a polynomial ring on $7-d$ variables.
\end{enumerate}
 \end{prop}
 \begin{pf}
 (1)
  It suffices to prove (1) for $L=\langle x_1,x_2 \rangle$ because $\Aut A_\RR$ acts transitively on 2-planes; however, after setting $x_1=x_2=0$, the relations $r_1$, $r_2$, and $r_3$, become the defining relations for 
  $\cA^{(+)}$,  and the relations  $r_4$, $r_5$, $r_6$, and $r_7$, imply that the 
 image of $x_3$ in $A/(x_1,x_2)$ commutes with the images of $x_4$, $x_5$, $x_6$, and $x_7$. 
 
 (2)
 Suppose $P$ is an associative 3-plane.   Because $\Aut A_\RR$ acts transitively on the set of
 associative 3-planes it suffices to prove (3) for $P=\langle x_1,x_2,x_3 \rangle$.  
However, it follows from the previous paragraph that $A/(x_1,x_2,x_3)$ is isomorphic 
 to   $\cA^{(+)}$. 
 
  Suppose $P$ is a non-associative 3-plane.
 Let $L$ be a 2-plane contained in $P$ and write $P=L \oplus \RR x$. By (1), $A/(L) \cong \cA^{(+)} \otimes k[t]$
 and the image of $x$ in   $\cA^{(+)} \otimes k[t]$ is a non-zero scalar multiple of $t+z$ where $z \in   \cA^{(+)}_1 -\{0\}$.
 Hence 
 $$
 {{A}\over{(P)}} \cong {{ \cA^{(+)}[t]}\over{(t+z)}} \cong  {{\cA^{(+)}[t]}\over{([t,y], [z,y] \; | \; y \in \cA^{(+)}_1)}}.
 $$
 However, by Lemma \ref{lem.CDV},  $\cA^{(+)}$ modulo the ideal  $([z,y] \; | \; y \in \cA^{(+)}_1)$ is a polynomial ring 
 on 4 variables so $A/(P)$ is too.  
 
 (3)
 If $\dim H \ge 4$, then $H$ contains a non-associative 4-plane so the result follows from (3). 
 \end{pf}

   \begin{prop}
  \label{prop.pdim1.1}
  Let $A=A_\RR$. If $L$ is a 6-plane   in $A_1$, then
  $$
  \pdim_A {{A}\over{ALA}} =2.
  $$
  \end{prop}
  \begin{pf}
 Since  $\Aut A_\RR$ acts transitively on lines in $V$, it acts transitively on their orthogonals. 
    We can therefore assume that $L=kx_2+\cdots + kx_7$.  
If $2 \le i \le 7$, then $[x_1,x_i] \in AL$ so $x_ix_1 \in AL$. Hence $AL=ALA$. 
Let $\mu:A \otimes L \to AL$ be the multiplication map. Since $A/AL$ is a polynomial ring in one variable,
the Hilbert series of $\ker \mu$ is
$$
6tH_A(t) - H_A(t) +(1-t)^{-1} = t^2H_A(t).
$$
But $\ker \mu$ contains $x_2 \otimes x_3-x_3 \otimes x_2 +x_4 \otimes x_5-x_5 \otimes x_4 +x_6 \otimes x_7-x_7 \otimes x_6$ and because $A$ is a domain the submodule of $\ker \mu$ generated by this element is 
has Hilbert series $t^2H_A(t)$. Hence $\ker \mu$ is the generated by this element, and the minimal projective
resolution of $A/ALA$ is 
$$
0 \to A(-2) \to A \otimes L \to A \to A/AL \to 0.
$$
It follows that $\pdim_A(A/AL)=2$.
\end{pf}

  \begin{prop}
  \label{prop.pdim1.1}
  Let $A=A_\RR$. If $L$ is a 5-plane   in $A_1$, then
  $$
  \pdim_A {{A}\over{ALA}} =1.
  $$
  \end{prop}
  \begin{pf}
 Since  $\Aut A_\RR$ acts transitively on the 2-planes in $V$, it acts transitively on their orthogonals. 
    We can therefore assume that $L=kx_3+\cdots + kx_7$. Let $L'=kx_4+\cdots + kx_7$.
   We will show that the multiplication map
   $$
  \mu: A \otimes_k(L + L'x_2+L'x_2^2+ \cdots ) \longrightarrow ALA
   $$
   is an isomorphism thus showing that $ALA$ is a free left $A$-module and hence that $A/ALA$ has 
   projective dimension one as required.
   
   First, to show that  $\mu$ is surjective we will show that the left ideal
   $$
   I:=A L + AL'x_2+AL'x_2^2+ \cdots,
$$
  which contains $AL$ and is contained in $ALA$, is a two-sided ideal and hence equal to $ALA$.
  Certainly $I$ is closed under right multiplication by $x_3,\ldots,x_7$ because those elements belong to $L$.
The element $ x_3x_2=x_2x_3+[x_4,x_5]+[x_6,x_7]$ belongs to $AL$
  so $Lx_2 \subset AL+L'x_2$; it follows that $Ix_2 \subset I$. Certainly 
  $x_3x_1=x_1x_3-[x_4,x_6]-[x_7,x_5]$   belongs to $AL$. The relations $r_4$, $r_5$, $r_6$, and $r_7$, show that $x_5x_1$, $x_4x_1$, $x_7x_1$, and $x_6x_1$ belong to $AL+L'x_2$. Hence $Lx_1 \subset AL+L'x_2 \subset I$. We will now show by induction on $n$ that $L'x_2^nx_1 \subset AL+AL'x_2+\cdots + AL'x_2^{n+1}$. We have just checked this is true for $n=0$ so suppose $n \ge 1$. Then
  $$
  L'x_2^nx_1 =L'x_2^{n-1}\big(x_1x_2 +[x_6,x_5]+[x_7,x_4]\big)
  $$
  which is contained in $\big( AL+AL'x_2+\cdots + AL'x_2^{n}\big)x_2+AL$. The induction argument therefore proceeds and it follows that $Ix_1 \subset I$.
  Thus $I=ALA$ and $\mu$ is surjective. 
  
  Let $D=L + L'x_2+L'x_2^2+ \cdots $. 
  To show that $\mu$ is injective it suffices to show that the Hilbert series of $ALA$ is equal to $H_A(t)H_D(t)$.
  Since $L'x_2^n \subset A_{n+1}$, $D=L \oplus L'x_2 \oplus L'x_2^2 \oplus \cdots$. Since $A$ is a domain
  $\dim_k(L'x_2^n)=\dim_kL'=4$. Hence $H_D(t)=5t+4t^2+4t^3+ \cdots = t+4t(1-t)^{-1} = (5t-t^2)(1-t)^{-1}$.
  
  Since $A/ALA$ is a polynomial ring in two variables,
  $$
  H_{ALA}(t)=H_A(t) - (1-t)^{-2} = (5t-t^2)(1-t)^{-1}H_A(t) = H_D(t)H_A(t)
  $$
  as required. This completes the proof.
  \end{pf}

     \begin{prop}
     [Piontkovski]
     \label{prop.pdim1.2}
  Let $A=A_\RR$. If $L$ is a co-associative 4-plane in $A_1$, then
  $$
  \pdim_A {{A}\over{ALA}} =1.
  $$
  \end{prop}
  \begin{pf}
Since $\Aut A_\RR$ acts transitively on co-associative 4-planes it suffices to prove the result for $L=kx_2+\ldots+kx_5$. That is what we will do. We will use a similar strategy to that in the proof of Proposition 
\ref{prop.pdim1.1}.

Let $J=AXA$. 
 Then $A/J$  is a commutative polynomial ring in three variables so has Hilbert series $(1-t)^{-3}$. The Hilbert series of $J$ is therefore 
\begin{align*}
H_J(t)= & H_A(t)- (1-t)^{-3}
\\
= &H_A(t) {{1}\over{(1-t)^3}}(1-3t+3t^2-t^3-1+7t-7t^2+t^3)
\\
= & H_A(t) {{4t}\over{(1-t)^2}}.
\end{align*}
Hence, to prove $J$ is free as a left $A$-module it suffices to show it is generated as a left $A$-module by a 
graded subspace whose Hilbert series is ${{4t}\over{(1-t)^2}}$.

Let $D$ be the linear span of the elements $\{x_1^ix_6^j \; | \; i,j \ge 0\}$.
Let $XD$ denote the linear span of $\{xd \; | x \in X, \, d \in D\}$. By Lemma \ref{lem.basis}, the multiplication map $X \otimes D \to XD$  is an isomorphism of graded vector spaces so $XD$ has Hilbert series 
$4t(1-t)^{-2}$.

Let $R$ be the subalgebra of $A$ generated by $x_1$ and $x_6$. 

The remainder of the proof is devoted to showing that 
$$
AXD = AXR =AXA = J.
$$
This will imply that the multiplication map $\phi:A \otimes X \otimes D \to J$ is a surjective homomorphism  between left $A$-modules having the same Hilbert series, and hence an isomorphism, 
thus showing that $J$ is a free left $A$-module.  
 
Certainly, $D$ contains $x_1D$ and $Dx_6$. Because $1 \in D$, $AX \subset AXD$.

Suppose $x_6^mx_1 \in D+AXD$. 
There is a relation of the form $x_6x_1=x_1x_6+r$ with $r \in AX$, so
\begin{align*}
x_6^{m+1}x_1  & =x^m_6(x_1x_6+r)
\\
&   \in (x_6^mx_1)x_6+ AX
\\
& \subset (D+AXD)x_6+AX
\\
& \subset D+AXD.
\end{align*}
Since $x_1 \in D+AXD$, it follows by induction on $m$ that $x_6^mx_1 \in D+AXD$ for all $m \ge 0$. Therefore
$x_1^\ell x_6^mx_1 \in  D+AXD$ for all $\ell, m \ge 0$, i.e., $$Dx_1 \subset D+AXD.$$
 
The next step is to show that $R_n \subset D+AXD$. This is certainly true for $n=0$ and if true for $n$,
then  
\begin{align*}
R_{n+1} & =R_n x_1 + R_n x_6
\\
&
 \subset (D+AXD)x_1+(D +AXD)x_6
 \\
 & \subset D+AXD.
\end{align*}
Therefore
$
R \subset D+AXD.
$
It follows that $AXD =AXR$.

Let $B$ be the subalgebra of $A$ generated by $x_1,\ldots,x_6$.

Since  $AXR(kx_1+\cdots +kx_6)  \subset AX+AXR  \subset AXR$ it follows that $AXB \subset AXR \subset AXB$, and hence $AXD=AXR=AXB$.
Therefore, to show that $AXD$ is equal to $AXA$, we need only show that $AXBx_7 \subset AXB$. Since $[x_7,-]$  is a derivation of $B$, 
\begin{align*}
AXBx_7  &  \subset AXB+ AXx_7B
\end{align*}
 But
$$
 \begin{array}{ccc}
 \,  x_2x_7 =x_7x_2 -  [x_1,x_4]  - [x_3,x_6]   &  \quad  &   x_4x_7 =x_7x_4 - [x_2,x_1]  -   [x_5,x_6]  
  \\
  \,  x_3x_7 =x_7x_3 - [x_5,x_1] - [x_6,x_2]   &  \quad &   x_5x_7 =x_7x_5 - [x_1,x_3] -  [x_6,x_4]
\end{array} 
$$
so $Xx_7 \subset AX +AXB = AXB$. Hence $AXBx_7 \subset AXB$. It follows that $AXB$ is a two-sided ideal
of $A$ and therefore $J=AXA=AXB=AXD$.
This completes the proof.
  \end{pf}

   \begin{prop}
   Suppose $k=\RR$. Let $L$ be a subspace of $A_1$ having dimension $\le 5$.
   Then the multiplication   map $A \otimes_k L \to AL$ is an isomorphism. In particular, 
   $AL$ is a free left $A$-module.
   \end{prop}
   \begin{pf}
   It is enough to prove this when $\dim_k L=5$ so we assume that is the case.
Since  $\Aut A_\RR$ acts transitively on the 2-planes in $V$, it acts transitively on their orthogonals. 
    We can therefore assume that $L=kx_2+\cdots + kx_6$. 
    Let $B=k[x_1,\ldots,x_6]$. Since $A$ is a free right $B$-module, $A \otimes_B I \cong AI$ for every 
    left ideal $I$ in $B$. It therefore suffices to show that  the multiplication   
    map $\mu:B \otimes_k L \to BL$ is an isomorphism of left $B$-modules.
    
Let $\cU$ be the set of words in $x_1,\ldots,x_6$ that do not contain $x_6x_1$ as a subword. Then $\cU$
is a basis for $B$. In particular, $\cU x_2 \sqcup \cdots \sqcup \cU x_6$ is linearly independent so
 $\mu$ is injective. But $\mu$ is surjective and hence an isomorphism. This completes the proof.
   \end{pf}

 \begin{prop}
 \label{prop.Hseries}
 Suppose $k=\RR$. If $x \in A_1-\{0\}$, then the Hilbert series of $A/(x)$ is $(1-t)^{-1} (1-5t+2t^2)^{-1}$.
 \end{prop}
 \begin{pf}
Since $\Aut A_\RR$ acts transitively on the set of lines in $V$, it suffices to prove the result for $x=x_7$.

 We first consider the algebra
 $$
R:={{\RR\langle x_1,x_2,x_3,x_4,x_5\rangle} \over{( [x_2,x_3] + [x_4,x_5] , [x_2,x_4] + [x_5,x_3]  )}} \,.
$$
 For brevity we will write $ij$ for $[x_i,x_j]$, $[ij,k]$ for $[[x_i,x_j],x_k]$, and
 $[i,jk]$ for $[x_i,[x_j,x_k]]$.  
 
The linear map $\d:R_1 \to R_2$ given by
 $$
 \begin{array}{lllll}
 \d(x_1)= [x_5,x_2]+[x_4,x_3]    & \quad &  \d(x_3)=[x_1,x_4]  & \quad &  \d(x_5)= [x_2,x_1]   \\
  \d(x_2)=[x_1,x_5]   & \quad &  \d(x_4)= [x_3,x_1]  & \quad &\phantom{ \d(x_5)= [x_2,x_1] }\\
\end{array} 
$$
extends to a derivation of $R$ because
 \begin{align*}
 \d\big( [x_2,x_3] + [x_4,x_5]  \big)  =   & \; [15,3]+[2,14]+[31,5]  +[4,21] \\
  =   & \; [35,1]+[1,24]  \\
  =   & \; 0 
\end{align*}
and
 \begin{align*}
 \d\big( [x_2,x_4] + [x_5,x_3]  \big)  = & \; [15,4]+[2,31]+[21,3]  +[5,14] \\
  =   & \; [45,1]+[1,32]  \\
  =   & \; 0.
\end{align*}
The Ore extension $R[x_6;\d]$ is therefore the quotient of the free algebra 
$\RR\langle x_1,\cdots,x_6\rangle$ modulo the relations 
$$
\begin{array}{lll}
 \,  \phantom{ .}\! [x_2,x_3] + [x_4,x_5] =0   &  \;  &  [x_2,x_4] + [x_5,x_3] =0   \\
  \,   \phantom{ ..  [x_6,x_4]+}  [x_6,x_1]= [x_5,x_2]+[x_4,x_3]   &  \;  &  
       \phantom{ ..  [x_6,x_4]+}  [x_6,x_2]= [x_1,x_5]     \\
 \,  \phantom{  .. [x_6,x_4]+}   [x_6,x_3]= [x_1,x_4]    &  \;  &     
     \phantom{ ..  [x_6,x_4]+}   [x_6,x_4]= [x_3,x_1]    \\
     \phantom{ ;.  [x_6,x_4]+}     [x_6,x_5]= [x_2,x_1].   &  \;  &  \phantom{ [x_2,x_4] + [x_5,x_3] =0  } \\
\end{array} 
$$
However, these 7 relations are defining relations for $A/(x_7)$, so $A/(x_7) \cong  R[x_6;\d]$.
The Hilbert series of  $A/(x_7)$ is therefore $(1-t)^{-1}H_R(t)$. 

A basis for $R$ is given by all words in $x_1,\ldots,x_5$ that do not contain $x_5x_4$ or 
$x_5x_3$ as a subword because
the relations for $R$ are  
\begin{align*}
x_5x_4 & =x_4x_5+ x_3x_2-x_2x_3 \\
x_5x_3 & =x_3x_5+ x_4x_2-x_2x_4.
\end{align*}
A word in $x_1,\ldots,x_5$  is {\sf good} if it does not contain $x_5x_4$ or $x_5x_3$ as a subword.
Let's write
$$
\begin{array}{ll}
U_n :=  \{\hbox{good words of length $n$ ending in $x_5$}\}   & \qquad u_n:=|U_n| 
\\
V_n :=  \{\hbox{good words of length $n$ not ending in $x_5$}\}   & \qquad v_n:=|V_n| 
\end{array}
$$
and 
$$
U(t)= \sum_{n=0}^\infty u_nt^n
\qquad \hbox{and} \qquad 
V(t)= \sum_{n=0}^\infty v_nt^n.
$$
Thus $H_R(t)=U(t)+V(t)$. 

The set of good words of length $n+1$ is the disjoint union 
$$
U_n x_1 \sqcup U_n x_2 \sqcup U_n x_5 \sqcup V_n x_1 \sqcup  V_n x_2 \sqcup  V_n x_3 \sqcup  V_n x_4 \sqcup  V_n x_5. $$
Therefore 
\begin{align*}
U_{n+1} & =  U_n x_5 \sqcup  V_n x_5 \\
V_{n+1}  & = U_n x_1 \sqcup U_n x_2 \sqcup   V_n x_1 \sqcup  V_n x_2 \sqcup  V_n x_3 \sqcup  V_n x_4  
\end{align*}
and $u_{n+1}=u_n+v_n$ and $v_{n+1}=2u_n+4v_n$.
Since $u_0 =0$ and $v_0=1$, it follows that 
$$
{{1}\over{t}} U(t)=U(t)+V(t)
\qquad \hbox{and} \qquad 
{{V(t)-1}\over{t}}=2U(t)+4V(t).
$$
It follows that $U(t)=t(1-5t+2t^2)^{-1}$, $V(t)=(1-t)(1-5t+2t^2)^{-1}$ and
$H_R(t) = (1-5t+2t^2)^{-1}$. The result follows. 
 \end{pf}

 Because $Bx_7 \subset B+x_7B$, $Ax_7B$ is closed under right multiplication by $x_7$ and is therefore a two-sided ideal. In other words, $Ax_7A=Ax_7B$.

 \section{Ideals of $A_\CC$ generated by subspaces of $A_1$}
 
 In the next proof we use the notation and results in sections \ref{sect.relns/C} and \ref{sect.SO}.
  
 \begin{prop}
 \label{prop.quo.rings}
 Suppose $k$ contains a square root of $-1$ and $\fchar k \ne 2$. 
 \begin{enumerate}
  \item 
 Let $L$ be a non-isotropic 2-plane such that $k+L$ is a subalgebra of ${\mathbb O}_k$. Then 
 $A/(L)$ is isomorphic to   the free algebra $k\langle u_1,u_2,u_3,v_1,v_2 \rangle$ modulo the relations
 $$
\qquad [u_1,u_2]=[u_2,u_3]=[u_3,u_1]=[v_1,v_2]=[u_1+v_2,u_2-v_1]=0.
 $$
  \item{}
 Let $L$ be an isotropic 2-plane in $\Im {\mathbb O}_k$. If
 \begin{enumerate}
   \item 
  $L^2 = 0$, then $A/ALA$   is isomorphic to $k \langle t,u_2,u_3,v_1,v_2 \rangle$ modulo the five relations
\begin{align*}
 & [u_2,v_2] =[t,u_2]= [t,v_2]=0
  \\
&  [t,v_1]-[u_2,u_3]=[t,u_3]-[v_2,v_1]=0.
\end{align*}
  \item 
  $L^2 \ne 0$, then $A/ALA$  is isomorphic to  $k \langle t,u_1u_2,u_3,v_2 \rangle$ modulo 
  the seven relations
  \begin{align*}
&  [u_2,v_2]=[t,u_1]=[t,u_2]=[t,u_3]=0
\\
& [u_1,u_2]=[u_2,u_3]=[t,v_2]-[u_3,u_1]=0;
\end{align*} 
\end{enumerate}  
\item{}
$A/(u_1,u_2,v_3) \cong k[t,v_1,v_2]*k[u_3]$ where $ k[t,v_1,v_2]$ is a commutative polynomial ring.
\end{enumerate}
 \end{prop}
 \begin{pf}
 (1)
By Corollary \ref{cor.nd-2-planes}, $\Aut {\mathbb O}_k$ acts transitively on the set of non-isotropic 2-planes $L$ such that $k+L$ is a subalgebra of ${\mathbb O}_k$ so it suffices to prove the result for one such $L$, say $L=kt+kv_3$. The result now follows from Proposition \ref{prop.relns/C}.

 (2)
 By Lemma \ref{lem.2-planes},  $\Aut A_k$ acts transitively on the set of null 2-planes and on the set of non-null isotropic 2-planes so it suffices to prove particular cases of (a) and (b). For (a) we may take 
$L=ku_1+kv_3$ and for (b) we may take $L=kv_1+kv_3$. The result now follows from Proposition \ref{prop.relns/C}.
 \end{pf}

  \section{$A$ is graded coherent}

    \subsection{}
  The main result in this section, that $A$ is graded coherent, is due to Dmitiri Piontkovski and I am grateful 
  for his generosity  in allowing me to include it here. His proof consisted of proving Proposition \ref{prop.pdim1.2} and then applying the following result.

\begin{prop}
 [Piontkovski 2006]
 \label{prop.coherent.0}
 \cite[Prop. 3.2]{Pion1} 
 A ring $A$ is right coherent if it has a right noetherian quotient ring $A/J$ such that $J$ 
  is free as a left $A$-module.
\end{prop}

  A ring is {\sf left coherent} if all its finitely generated left ideals are finitely presented. 
A graded ring is {\sf graded left coherent} if the following equivalent conditions hold:
\begin{enumerate}
  \item 
  all its finitely generated graded left ideals are finitely presented;
  \item 
   every finitely generated graded submodule of a finitely presented graded module is 
  finitely presented;
  \item 
  the category of finitely presented graded left $R$-modules is abelian.
\end{enumerate}

Suppose $R$ is a connected graded ring. A finitely presented graded left $R$-module will be called a
{\sf coherent} module, and we write $\coh R$ for the abelian category of graded coherent left modules. 

Suppose $R$ is a connected graded $k$-algebra.  If the minimal projective resolution 
of a graded left $R$-module $M$ begins
$$
\cdots \to R \otimes_k V_1 \to R \otimes_k V_0 \to M \to 0,
$$
then $V_i \cong \Tor_i^R(k,M)$ so $M$ is coherent if and only if $\Tor_i^R(k,M)$ is finite 
dimensional for $i=0$ and $i=1$. It follows that an extension of coherent modules is coherent. If $R_{\ge 1}$ is a finitely generated left ideal, then the trivial module 
$k=R/R_{\ge 1}$ is coherent. It is also, up to twisting, the only
simple graded module, so it follows that every finite dimensional graded $R$-module is coherent.

 If a left coherent ring has finite left global homological dimension all its finitely 
 presented modules have a finite projective resolution by finitely presented projective modules.

 \begin{prop}
 [Piontkovski]
 \label{prop.coherent}
$A$ is graded coherent on the right and the left.
\end{prop}
\begin{pf}
Propositions \ref{prop.pdim1.1} and \ref{prop.pdim1.2} exhibit two-sided ideals $J$ such that $A/J$
is noetherian and $J$  is free as a left $A$-module. By Proposition \ref{prop.coherent.0}, $A$ is right coherent.
Since $A$ is an enveloping algebra it is isomorphic to its opposite ring, so it is also left coherent. 
\end{pf}

I do not know if $A \otimes_k A$ is graded coherent.

\section{Non-commutative geometry}
\label{sect.nag}

\subsection{}
We will prove a result for $A$ that is similar to a result for the projective plane. Let $\PP^2$ be the 
projective plane over a field $k$. Let $\sD^b(\coh \PP^2)$ be the bounded derived category of coherent
sheaves on $\PP^2$. The locally free sheaf $\cT:=\cO \oplus \cO(1) \oplus \cO(2)$ is a {\it tilting sheaf}
which means that it generates $\sD^b(\coh \PP^2)$ and that $\Ext^i(\cT,\cT)=0$ for $i \ne 0$. General
theory implies that the functor $\RHom(\cT,-)$ is an equivalence between $\sD^b(\coh \PP^2)$ and 
$\sD^b(\mod E)$ where $E$ is the endomorphism ring of $\cT$ and $\mod E$ is the category of finite dimensional $E$-modules. One can be explicit about $E$, namely
$$
E = \begin{pmatrix}
   S_0   & S_1 & S_2   \\
    0  &   S_0   & S_1 \\
    0 &  0  &   S_0
\end{pmatrix}
$$
where $S=k[X,Y,Z]$ is the polynomial ring and $S_i$ is its degree $i$ component. The algebra 
$E$ is also the path algebra of a quiver with relations, and the equivalence of categories may be stated
as an equivalence between $\sD^b(\coh \PP^2)$ and the bounded derived category of representations of that
quiver. Then quiver is
$$
 \UseComputerModernTips
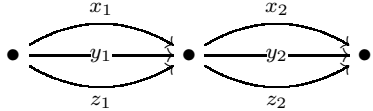
\xymatrix{
\bullet \ar@/^1pc/[rr]^{x_1}  \ar@/_1pc/[rr]_{z_1} \ar[rr] | {y_1} && \bullet \ar@/^1pc/[rr]^{x_2}  \ar[rr] | {y_2}  \ar@/_1pc/[rr]_{z_2}  && \bullet
}
$$
with  relations 
$$
y_2x_1-x_2y_1 = y_2z_1-z_2y_1 =z_2x_1-x_2z_1  = 0.
$$
The vertices represent $\cO$, $\cO(1)$, and $\cO(2)$, and the arrows represent bases for $\Hom(\cO,\cO(1))$ and $\Hom(\cO(1),\cO(2))$, and the relations arise from the composition map
$$
\Hom(\cO(1),\cO(2)) \otimes \Hom(\cO,\cO(1)) \to \Hom(\cO,\cO(2)).
$$

\subsection{}
The category of finite dimensional graded $A$-modules, $\fdim A$,  is a full subcategory of 
$\coh A$. It is a Serre subcategory so we may form the quotient category and, following Artin and Zhang \cite{AZ}, implicitly define a non-commutative space $X:=\Projnc A$ by declaring that the category of ``coherent sheaves'' on $X$ is
$$
\coh X:= {{\coh A}\over{\fdim A}}.
$$
Because $\gldim A=3$, $X$ has cohomological dimension two. An argument like that in \cite{MS}   
shows that 
$$
K_0(\coh X) \cong {{\ZZ[t,t^{-1}]}\over{(1-7t+7t^2-t^3)}}   \cong \ZZ^3.
$$

The next result is analogous to many other such theorems. The key point is that 
 $\cO$, $\cO(1)$, and $\cO(2)$, the images of $A$, $A(1)$, and $A(2)$, in $\coh X$, generate 
 $\sD^b(\coh X)$ and $\Hom_X(\cO(i),\cO(i-1))=0$.

\begin{thm}
\label{thm.equiv}
There are equivalences of categories 
$$
\sD^b(\coh X) \equiv \sD^b(\rep Q)  \equiv \sD^b(\mod E)
$$
where $Q$ is the quiver
$$
 \UseComputerModernTips
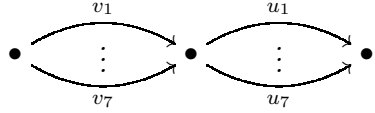
\xymatrix{
\bullet \ar@/^1pc/[rr]_{\vdots}^{v_1}  \ar@/_1pc/[rr]_{v_7} && \bullet \ar@/^1pc/[rr]_{\vdots}^{u_1}   \ar@/_1pc/[rr]_{u_7}  && \bullet
}
$$
with  relations 
\begin{equation}
\label{quiv.relns}
\begin{array}{cc}
u_2v_3-u_3v_2+u_4v_5-u_5v_4+u_6v_7-u_7v_6  & \, =\, 0, \\
u_3v_1-u_1v_3+u_4v_6-u_6v_4+u_7v_5-u_5v_7 & \, =\, 0, \\
u_1v_2-u_2v_1+u_7v_4-u_4v_7+u_6v_5-u_5v_6 & \, =\, 0, \\
u_5v_1-u_1v_5+u_6v_2-u_2v_6+u_3v_7-u_7v_3& \, =\, 0, \\
u_1v_4-u_4v_1+u_2v_7-u_7v_2+u_3v_6-u_6v_3 & \, =\, 0, \\
u_7v_1-u_7v_1+u_2v_4-u_4v_2+u_5v_3-u_3v_5 & \, =\, 0, \\
u_1v_6-u_6v_1+u_5v_2-u_2v_5+u_4v_3-u_3v_4 & \, =\, 0,
\end{array}
\end{equation}
and $E$ is its path algebra
\begin{align*}
\End_{\gr A} (A \oplus A(1) \oplus A(2)) 
& 
\cong \begin{pmatrix}
   A_0   & A_1 & A_2   \\
    0  &   A_0   & A_1 \\
    0 &  0  &   A_0
\end{pmatrix}
\\
& 
= \begin{pmatrix}
   k   & (\Im {\mathbb O})^* & \coker \mu^*   \\
    0  &   k   & (\Im {\mathbb O})^* \\
    0 &  0  &   k
\end{pmatrix}.
\end{align*}
The group ${\bf G}(k)$ of type $G_2$ acts as automorphisms of $E$.
\end{thm}
\begin{pf}
The proof follows a well-worn path. The category $\sD^b(\coh X)$ is generated by $\cT:=\cO \oplus
\cO(1) \oplus \cO(2)$ and $\Ext_{\coh X}(\cT,\cT)=0$ for $i > 0$; i.e., $\cT$ is a tilting sheaf. Hence 
$\RHom(\cT,-)$ is an equivalence from $\sD^b(\coh X)$ to $\sD^b(\mod \End \cT)$. 
The equivalence of categories follows from \cite[Thm. 4.14]{MM}. The other facts follow from earlier results in this paper. 
\end{pf}

\subsection{The $(1,1,1)$ moduli space for $Q$}

We now identify  $\PP^6$ with the space of lines in $\Im {\mathbb O}_k$ and write $[u]$ for the point in 
$\PP^6$ corresponding to the line $ku$ in $\Im {\mathbb O}_k$.

Up to isomorphism an  indecomposable representation of $Q$ having dimension vector 
$(1,1,1)$ determines, and is determined by, a point in $\PP^6 \times\PP^6$. 
Let  $\cM_{(1,1,1)}$ be the closed subvariety of $\PP^6 \times\PP^6$
consisting of points $((v_1,\ldots,v_7), (u_1,\ldots,u_7))$ satisfying the equations (\ref{quiv.relns}).

\begin{prop}
\label{prop.moduli}
$$
\cM_{(1,1,1)}=\{([u],[v]) \; | \; \Im(uv)=0\} \subset \PP^6 \times \PP^6.
$$ 
\end{prop}
\begin{pf}
Suppose $(v_1,\ldots,v_7), (u_1,\ldots,u_7) \in k^7-\{0\}$ and define
$$
u:= \sum_{i=1}^7v_io_i
\qquad v:= \sum_{i=1}^7 u_io_i.
$$
The quadratic terms in (\ref{quiv.relns}) are then the coefficients of $o_1,\ldots,o_7$ in the product $vu$,
so a point  $([u],[v]) \in \PP^6 \times \PP^6$ belongs to $\cM_{(1,1,1)}$ if and only if $\Im(vu)=0$ or, equivalently, $\Im(uv)=0$. 
\end{pf}

If ${\mathbb O}_k$ is a division algebra, $\cM_{(1,1,1)}$ is the diagonal copy of $\PP^6$ in $\PP^6$.
For example, over $\RR$, $\cM_{(1,1,1)} \cong \RR\PP^6 \cong S^6/\!\!\sim$ where we identify 
antipodal points on the 6-sphere.

\begin{prop}
Suppose $k$ is algebraically closed and of characteristic not 2.   Let $\pi_1$ and $\pi_2$ be the projections $\pi_1\big([u],[v]\big)=[u]$ and  $\pi_2\big([u],[v]\big)=[v]$, 
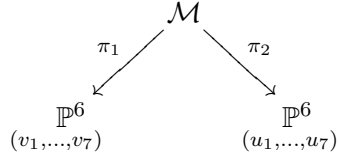
$$
 \UseComputerModernTips
\xymatrix{
& \cM \ar[dl]_{\pi_1} \ar[dr]^{\pi_2}
\\
\PP^6   \save []+<-0.2cm,-0.35cm>*\txt<5pc>{$\scriptstyle{(v_1,\ldots,v_7)}$} \restore &&
 \PP^6   \save []+<-0.1cm,-0.35cm>*\txt<5pc>{$\scriptstyle{(u_1,\ldots,u_7)}$} \restore
}
$$
Then $\pi_1^{-1}([u])= \{[u]\} \times \PP(E_u)$ where 
\begin{align*}
E_u := & \{v \in \Im {\mathbb O}_k \; | \; \Im(uv)=0\}
\\
= & \{v \in \Im {\mathbb O}_k \; | \; \phi(u,v,-)\equiv 0\}.
\end{align*}
Furthermore,
\begin{enumerate}
  \item 
 if $ku$ is not isotropic, then $\pi_1^{-1}([u])=([u],[u])$.
  \item 
If $ku$ is isotropic, then 
\begin{enumerate}
  \item 
  $E_u$ is an isotropic 3-plane and
  \item 
 $ E_u := \{v \in {\mathbb O}_k \; | \; uv=0\}$.
\end{enumerate}
\end{enumerate}
\end{prop}
\begin{pf}
(1)
In this case, $u$ is a unit so $E_u=ku^{-1}=ku$.

(2)
(a)
By Lemma \ref{prop.lines}, $\Aut {\mathbb O}_k$ acts transitively on the set of isotropic lines and
therefore  on the set of associated subspaces $E_u$. It therefore suffices to show that $\dim E_u=3$ for one particular choice of $u$. Let $u_1=o_2+io_3$. It follows from (\ref{mult.table.SO}) that
  $K_{u_1}=ku_1+kv_2+kv_3 \cong k^3.$  
  
(b)
Since $E_u$ is isotropic, $\Real(uv)=0$ for all $v \in E_u$. But $\Im(uv)=0$ by definition of $E_u$ so $uv=0$.  
\end{pf}

The locus in $\PP^6$ consisting of points $[u]$ such that $u^2=0$ is the smooth quadric 
$x_1^2+\cdots  +x_7^2=0$. If $k$ is algebraically closed and $\fchar k \ne 2$, then $\cM_{(1,1,1)}$ has two components, one the diagonal copy of $\PP^6$ in $\PP^6 \times \PP^6$, and the other a $\PP^2$-bundle over
the smooth quadric $x_1^2+\cdots  +x_7^2=0$. The dimension of the latter component is 7.

\subsubsection{} 
If $k$ is algebraically closed and of characteristic not 2, then every null 3-plane is of the form $E_u$ for some $u \in \Im {\mathbb O}_k$. To see this, suppose $L$ is a null 3-plane and let $ku$ be a line in $L$. Then $ku$
is isotropic and $uL=0$ so $L=E_u$.  

\subsubsection{}
If $k$ is algebraically closed and of characteristic not 2, then $\Aut {\mathbb O}_k$ acts transitively on the set of null 3-planes (we observed in the course of proving Proposition \ref{prop.moduli} that $\Aut {\mathbb O}_k$ acts transitively on the set of 3-planes of the form $E_u$). Hence if $L$ is a null 3-plane in $\Im {\mathbb O}_k$, then $A/ALA$ is isomorphic to $k \langle t,u_2,u_3,v_1 \rangle$ subject to the relations
$$
[t,u_2]=[t,u_3]=[t,v_1]-[u_2,u_3]=0.
$$

\subsubsection{}
We note that $\CC u_1+\CC u_2 +\CC u_3$ is an isotropic 3-plane in $\Im {\mathbb O}_\CC$ 
that is not null. Also, $A_\CC/(u_1,u_2,u_3)$ is the commutative polynomial ring on 4 variables.

\subsubsection{Point modules}

Our basic reference for point modules is \cite{ATV}. A point module over $A$ is a graded $A$-module 
that is generated as an $A$-module by its degree zero component and has Hilbert series $(1-t)^{-1}$.
If $M=ke_0 \oplus ke_1 \oplus \cdots$ is a point module for $A$, with $\deg e_i=i$, 
then there is for each $i$ a unique $\l_i \in A_1^*$ such that $x.e_i=\l_i(x) e_{i+1}$ for each $x \in A_1$. The 
image of $\l_i \in \PP(A_1^*)$ does not depend on the choice of $e_i$s so the point module determines, and is determined up to isomorphism, by the sequence $\l_0,\l_1,\ldots$ of points in $\PP(A_1^*)$. 

The points in $\cM_{(1,1,1)}$ correspond to truncated point modules of length 3 for $A$. 
If $u \in \Im {\mathbb O}_k -\{0\}$ and $u^2 \ne 0$, then $u$ corresponds to a point module that corresponds
to a point on $\Proj S(\Im  {\mathbb O}_k^*)=\PP(\Im {\mathbb O}_k)$. 
However, if $u \in \Im {\mathbb O}_k -\{0\}$ and $u^2 = 0$, there is a point module that corresponds to a 
sequence $([u],[v],[w],[x],\ldots)$ of points in $\PP(\Im {\mathbb O_k})$ such that $0=uv=vw=wx=\cdots $ and each such sequence determines a point module for $A$.

 \section{Special elements in $A$}
 
 The next lemma implies that over $\CC$ there are no ``non-trivial'' relations of the form $ab=cd$
 between elements $a,b,c,d \in A_1$; i.e., $R$ contains no rank-two tensors.
 
We point this out because relations of the form $ab=cd$ between degree-one elements  play an important role in analyzing the homogeneous coordinate rings of various non-commutative analogues of $\PP^2$ and 
$\PP^3$. For example, if $S$ is a 3- or 4-dimensional Sklyanin algebra and $ab=cd$ is a non-trivial relation between elements of $S_1$, then $S/Sb+Sd$ is a 
 point module (or line module), and all point modules (line modules) are of this form.   
  
 \begin{lem}
 \label{lem.no.rk2} 
Let $R \subset \Im {\mathbb O}_k^* \otimes \Im {\mathbb O}_k^*$ be the space of relations and 
$W$ the element defined in (\ref{defn.W}). 
 \begin{enumerate}
  \item 
  The rank of $W$ is 7.
    \item{}
  If ${\mathbb O}_k$ is a division algebra, for example $k=\RR$, 
  every  non-zero element in $R$ has rank six.
  \item 
  If $k$ is algebraically closed of characteristic $\ne 2$, 
  the rank of every non-zero element in $R$ is $\ge 4$ and equality can occur.
\end{enumerate}
\end{lem}
\begin{pf}
Let $V=\Im {\mathbb O}_k^*$. 

(1)
Under the isomorphism $V^{\otimes 3} \to \Hom(V^*\otimes V^*,V)$ we consider $W$
as a linear map $V^*\otimes V^* \to V$ by
$$
W(\l,\mu):= \sum_{ijk}\ve^{ijk} x_i\l(x_j)\mu(x_k).
$$
Each $x_i$, $1 \le i \le 7$, is in the
image of $W$ because
$$
W(o_j, o_k) = \pm o_i
$$
whenever $ijk$ is a line in the Fano plane. 
Hence $\rank W=7$.

(2)
The hypothesis implies  there are no isotropic lines in $\Im {\mathbb O}_k$ and therefore
 $\Aut {\mathbb O}_k$ acts transitively on the lines in $\Im {\mathbb O}$ and hence on the lines in $R=\im \mu^*$. But one relation has rank 6 so all non-zero relations have rank 6.

(3)
Since $k$ is algebraically closed it contains a square root of $-1$. In that case, section \ref{sect.relns/C}  exhibits several rank 4 elements in $R$.

Suppose $R$ has a non-zero element of rank $< 4$.
Since $R$ consists of skew-symmetric tensors it contains a non-zero element of the form
 $u \otimes v - v \otimes u$. The isotropic lines  in $\Im {\mathbb O}_k$ are the points of a smooth quadric hypersurface in $\PP^6$. This quadric is a union of projective 3-planes but does not contain a 
 projective 4-plane. 
The 5-plane $(k u + k v)^\perp$ therefore contains a non-isotropic line, $\CC w$ say. But 
$\Aut {\mathbb O}_k$ acts transitively on the set of non-isotropic lines (Proposition \ref{prop.lines})
 so $w=g\cdot x_7$ for some $g \in \Aut {\mathbb O}_k$.
Hence $u$ and $v$ belong to $(g\cdot x_7)^\perp=g \cdot (k x_1+\cdots +k x_6)$. 
The subalgebra of $A$ generated by $(g\cdot x_7)^\perp$ is therefore  isomorphic to 
the algebra $B$ in (\ref{eq.B.alg}). However, the hypothesis
 that $u \otimes v-v \otimes u$ is a relation says that $g\cdot (u \otimes v - v \otimes u)$ is a non-zero
 multiple of the relation $[x_1,x_6]+[x_5,x_2]+[x_4,x_3]$ for $B$ and therefore of rank 6. This is absurd so we conclude $R$ does not contain a non-zero relation of rank $<4$. .
 \end{pf}

  An element $z$ is {\sf normal} in $A$ if $zA=Az$.

\begin{prop}
\label{prop.center}
Let $k$ be any field of characteristic zero.   
\begin{enumerate}
  \item 
  The only normal elements in $A$ are the elements in $k$.
  \item 
  The center of $A$ is equal to the base field, $k=A_0$. 
\end{enumerate}
\end{prop}
\begin{pf}
Let $\kol$ be an algebraic closure of $k$. It suffices to prove the result for $A \otimes_k \kol$ so we can, and do, assume that $k$ is algebraically closed. 
It suffices to show there are no homogeneous normal elements of degree $\ge 1$.

As before, let $B=k[x_1,\ldots,x_6]$. Then $A=B[x_7;\d]$. 

Write $x=x_7$. Each element $a$ in $A$ can be written as  $a=\sum_{i=1}^n b_i x^i$
for unique elements $b_i$ in $B$ with $b_n \ne 0$. 
We say that such an $a$ has $x$-degree $n$, and write 
$\deg_x(a):=n$. 

Suppose $b_1, \ldots ,b_n \in B$ are such that
\begin{equation}
\label{central.elt?}
z:= \sum_{i=1}^n b_i x^i
\end{equation}
is a homogeneous normal element. Without loss of generality, we assume that $b_n\ne 0$.

Let $c$ be a homogeneous element in $B$. Then there is a unique element $\widehat{c} \in A$ such that
 $cz=z\widehat{c}$. 
Since $A$ is an Ore extension of $B$ by $x$, $\deg_x(uv)=\deg_x(u)+\deg_x(v)$ for 
every pair of elements $u,v \in A$. Hence $\deg_x(\widehat{c})=\deg_x(b)=0$, which says that 
$\widehat{c}$ belongs to $B$. Therefore $Bz \subset zB$. In fact, since $z$ is homogeneous and $A$ is a 
domain, $B_jz=zB_j$ for all $j$.

Now
\begin{align*}
z\widehat{c} & = b_nx^n\widehat{c}+b_{n-1}x^{n-1}\widehat{c} + \hbox{terms of lower $x$-degree }
\\
& = b_n\widehat{c}x^n+ \big[b_{n-1}\widehat{c}  + n b_n \d(\widehat{c})\big]x^{n-1} + \hbox{terms of lower $x$-degree}
\end{align*}
so  $cb_n=b_n \widehat{c}$ and $cb_{n-1} = b_{n-1}\widehat{c} + n b_n \d(\widehat{c})$. Because $c$ was arbitrary, the fact that 
$cb_n=b_n \widehat{c}$ implies that $Bb_{n}\subset b_{n}B$. Since $B$ is a domain, $Bb_{n}$
has the same Hilbert series as $b_{n}B$ so $Bb_{n}  =  b_{n}B$.  But
 the only normal elements in $B$ are the elements in $k$ (\cite[Thm. 0.3(3)]{Z}) so 
$b_n \in k$.  We may therefore assume that $b_n=1$.
Since $z$ is homogeneous it follows that $b_i \in B_{n-i}$ for all $i$. In particular, $b_{n-1} \in B_1$. 

The equality  $cb_{n-1} = b_{n-1}\widehat{c} + n b_n \d(\widehat{c})$ now becomes $cb_{n-1} = b_{n-1}\widehat{c} + n  \d(\widehat{c})$.
But $b_{n-1}z=z\widehat{b_{n-1}}$ for some $\widehat{b_{n-1}} \in A$ so $b_{n-1}(b_{n-1} - \widehat{b_{n-1}} )= n \d(\widehat{b_{n-1}})$.
More explicitly, 
\begin{equation}
\label{elt.rank.4}
b_{n-1} \otimes (b_{n-1} - \widehat{b_{n-1}}) -  n (x_7 \otimes\widehat{b_{n-1}} -\widehat{b_{n-1}}\otimes x_7) \in R.
\end{equation}
But every non-zero element of $R$ has rank $\ge 4$ by Lemma 
\ref{lem.no.rk2} so the element in (\ref{elt.rank.4}) is zero. However, the three terms belong to 
$B_1 \otimes B_1$, $x_7 \otimes B_1$ and $B_1 \otimes x_7$, so each individual term is zero. Since 
$n \ne 0$ it follows that $\widehat{b_{n-1}}=0$ and hence $b_{n-1}=0$.

Since  $B_1 z=zB_1$, every element in  $B_1$ is of the form $\widehat{c}$ for some $c \in B_1$; i.e., 
$cz=z\widehat{c}$.  
 Since $b_{n-1}=0$ we have
 \begin{align*}
z\widehat{c} & = b_nx^n\widehat{c} + \hbox{terms of $x$-degree $\le n-2$}
\\
& = b_n\widehat{c}x^n+   n b_n \d(\widehat{c})x^{n-1} + \hbox{terms of lower $x$-degree}
\end{align*}
 But $n \ne 0$, so  $\d(\widehat{c})=0$. But $\d=[x_7,-]$ has no kernel in $B_1$ so we conclude that  
 no such $z$ can exist.
\end{pf}

\subsection{A special degree 4 element}

By the Poincar\'e-Birkhoff-Witt theorem applied to $A$, the element 
\begin{equation}
\label{Q.A4}
Q:= \sum_{1 \le p<q \le 7} [x_p,x_q]^2 \, \in \, A_4 
\end{equation}
is non-zero. This element has some interesting properties when the base field is $\RR$. For example, 
it is $G_2$-invariant and belongs to {\it every} subalgebra of $A_\RR$ generated by a codimension-one subspace of $A_1$ (Proposition \ref{prop.Q.everywhere}). 

$Q$ also appears in the identity (\ref{main.eq}) which is the key step in the proof of the next result. Richard Eager informed me that both the identity (\ref{main.eq}) and the result itself was already known to physicists---see
 \cite[Eq. (17)]{FST}.   However, they did not give a proof so we do that here.

  \begin{prop}
 \label{prop.mods}
 Let $M$ be a finite dimensional left $A_\RR$-module. If $M$ is isomorphic to its dual as a module over the Hopf algebra $A_\RR$, then $M$ is annihilated by $[x_i,x_j]$ for all $i$ and $j$. In other words, $M$ is actually a module over the commutative ring $S(V)$. 
 \end{prop}
 \begin{pf}
 The hypothesis implies that $M$ has a basis with respect to which each element in $A_1$ acts on
 $M$ as  a skew symmetric matrix.
 
 We consider $M$ as a module over the free algebra $TV$. 
We write $\Tr_M(a)$ for the trace of an element $a \in TV$  acting on   $M$. 
 
Since $M$ is annihilated by the relations $r_1,\ldots,r_7$ it is annihilated by   
$$
 \sum_{i=1}^7  r_i^2  =  \sum_{1 \le p<q \le 7} [x_p,x_q]^2 +d = Q+d
 $$
 where $d$ is the sum of all terms of the form $[x_p,x_q][x_r,x_s]$ where $\{p,q,r,s\}$ range over all 
 4-tuples of distinct points of the Fano plane, no three of which are colinear. 
 
 We will now show $\Tr_M(d)=0$. Fix  distinct points $\{p,q,r,s\}$ in the 
 Fano plane, no three of which are colinear, and let $d_{pqrs}$ be the part of $d$ consisting of the
 words in which the four distinct letters $\{x_p,x_q,x_r,x_s\}$ appear.  
 If we write $[pq][rs]$ for the product of commutators $[x_p,x_q][x_r,x_s]$, then
 $$
 d_{pqrs}=[pq][rs] +[rs] [pq]+[pr][qs] +[qs][pr] +[ps][qr] +[qr][ps].
 $$
An explicit calculation shows that $\Tr_M(d_{pqrs})=0$; this can also be proved by analyzing the cofficients of the form $\ve^{ipr}\ve^{iqs} + \ve^{jpq}\ve^{jrs}+\ve^{krq}\ve^{kps}$ that appear in $\sum_{i=1}^7 r_i^2$.

Since $d$ is the sum of all $d_{pqrs}$ as  $\{p,q,r,s\}$ runs over all non-colinear 4-tuples in the Fano plane it follows that $\Tr_M(d)=0$. 
Therefore 
 \begin{equation}
\label{main.eq}
 \sum_{1 \le p<q \le 7} \Tr_M([x_p,x_q]^2)  =  \sum_{i=1}^7 \Tr_M(r_i^2)  = 0.
\end{equation}

Now fix a basis for $M$ and consider the matrix form of the elements of $A$ acting on $M$ 
with respect to that basis. 

Since $x_p$ and $x_q$ are skew-symmetric, so is $[x_p,x_q]$. Hence  the  eigenvalues of 
$[x_p,x_q]$, as for any real skew-symmetric matrix, are either zero or purely imaginary. The eigenvalues of $[x_p,x_q]^2$ are therefore real and $\le 0$. It now follows from (\ref{main.eq}) that the only eigenvalue of $[x_p,x_q]^2$, and hence the only eigenvalue of $[x_p,x_q]$, on $M$ is  zero. 
Hence $[x_p,x_q]$ acts as zero on $M$.
  \end{pf}
  
  \begin{lem}
 Let $k=\RR$. Then $Q$ is fixed by $G^c_2$.
 \end{lem}
 \begin{pf}
 Since $G^c_2$ is connected it suffices to show that the Lie algebra of $G^c_2$ annihilates $Q$. 
 The Lie algebra of $G_2^c$ acts as degree-preserving derivations of $A$ so it suffices to show that every 
 such derivation annihilates $Q$. We do this is section \ref{sect.Der.A} below.
 \end{pf}

 \begin{prop}
 \label{prop.Q.everywhere}
 Let $k=\RR$.
 Let $U$ be a codimension one subspace of $A_1$. 
 Then the subalgebra of $A_\RR$ generated by $U$ contains the element $Q$.
\end{prop}
\begin{pf}
Since $G_2$ acts transitively on the lines in $\Im({\mathbb O})$  it acts transitively on the codimension-one subspaces of $A_1$.  Since $Q$ is a $G_2$-invariant it suffices to show that $Q$ is contained in the subalgebra $B=\RR[x_1,\ldots,x_6]$ of $A$. If $b \in B$, then $[b,x_7] \in B$ by Proposition \ref{prop.Ore}. Hence 
$$
 \sum_{1 \le p \le 6} [x_p,x_7]^2 \in B.
$$
Since $B$ also contains
$$
 \sum_{1 \le p<q \le 6} [x_p,x_q]^2 
$$
it contains $Q$.
\end{pf}

\subsection{Homogeneous derivations of $A$}
\label{sect.Der.A}

Because $\Aut_{\sf gr} A$ contains an algebraic group of type $G_2$ the corresponding
Lie algebra of type $G_2$ will act as degree-preserving derivations of $A$.  
Let $$\Der_{\sf gr}A_k$$ denote the space of degree-preserving $k$-linear derivations of $A$.

\begin{prop}
\label{prop.ders}
Let $j$ and $k$ be the endpoints of an arrow $j \to k$  in the Fano plane. (There are 21 such pairs $jk$.)
\begin{enumerate}
\item
There is a unique labelling of the points in the Fano plane such that 
 $ijk$, $ipq$, $irs$, and $jpr$,  are directed lines.
  \item 
There is an element $\d_{jk}$ in $\Der_{\sf gr}A_k$ acting as follows:
$$
\begin{array}{|r|c|c|c|c|c|c|c|}
\hline
x= &x_i & x_j & x_k & x_p & x_q & x_r & x_s
\\
\hline
\d_{jk}(x)= & 0 & 0 & 0 & x_q & -x_p & -x_s & x_r \\
\hline
\end{array}
$$
\item{}
The derivations $\d_{jk}$ span $\Der_{\sf gr}A_k$. 
\item{}
These derivations satisfy the identity
$$
\d_{jk}+\d_{pq}+\d_{rs}=0.
$$
\item{}
$\Der_{\sf gr}A_k$ is a Lie algebra   of type $G_2$.
\item
$\d_{ij}$, $\d_{jk}$, and $\d_{ki}$, span a Lie subalgebra 
isomorphic to $\fso(3)$ with explicit relations
$$
[\d_{ij},\d_{jk}]=2\d_{ki}, \quad [\d_{jk},\d_{ki}]=2\d_{ij}, \quad [\d_{ki},\d_{ij}]=2\d_{jk}.
$$
\end{enumerate}
\end{prop}

The following table lists the 21 derivations.  The labelling $\d_{jk}$ is as indicated in Proposition 
\ref{prop.ders}. It is the subscripts of the $x_i$s that carry information so we only 
enter the subscripts into the table with the understanding that an entry $i$ stands for $x_i$ and an entry $\iol$ stands for $-x_i$.
Furthermore, if $\d(x_i)=0$ we  leave the corresponding entry in the table blank.
 $$
\begin{array}{|c|c|c|c|c|c|c|c|c|c|c|}
\hline
ijk &ipq& irs &\d_{jk}&x_1 & x_2 & x_3 & x_4 & x_5 & x_6 & x_7
\\
\hline
123 & 145 & 167  & \d_{23} &   &   &   & 5 & \4ol & \7ol & 6 \\
145 & 167 & 123  & \d_{45} &   & \3ol & 2 &   &   & 7 & \6ol \\
167 & 123 & 145  & \d_{67} &   &  3 & \ol2  & \5ol & 4 & &  \\
231 & 275 & 246  & \d_{31} &   &   &   & \6ol & \7ol & 4 & 5 \\
246 & 231 & 275  & \d_{46} &  \3ol &   & 1  & &  7 &   & \5ol \\
275 & 246 & 231  & \d_{75} & 3 &   & \1ol  & 6 &   &  \4ol & \\
312 & 365 & 374  & \d_{12} &   &   &   & 7 & \6ol & 5 & \4ol \\
365 & 374 & 312  & \d_{65} & \ol2  & 1  &   & \7ol&   &   &4 \\
374 & 312 & 365  & \d_{74} & 2& \1ol &   & & 6  & \5ol  & \\
451 & 437 & 462  & \d_{51} &   & 6 & 7 &   &   & \ol2 & \3ol \\
437 & 462 & 451  & \d_{37} & 5 & \6ol  &   & & \1ol  & 2  & \\
462 & 451 & 437  & \d_{62} & \5ol  &   &  \7ol & & 1  &   & 3\\
514 & 527 &  536 & \d_{14} &   &  7 & \6ol  &  &  & 3 & \ol2 \\
527 & 536 & 514  & \d_{27} & \4ol  &   &  6 & 1&   & \3ol  & \\
536 & 514 & 527  & \d_{36} & 4  & \7ol  &   & \1ol&   &   &2 \\
671 & 653 & 624  & \d_{71} &   & \4ol  & \5ol  &2 & 3 &   & \\
653 & 624 & 671  & \d_{53} & 7 &  4 &   & \ol2&   &   & \1ol \\
624 & 671 & 653  & \d_{24} &  \7ol &   & 5  & &\3ol   &   & 1\\
716 & 743 & 752  & \d_{16} &   & 5  & \4ol  &3 & \ol2  &   & \\
752 & 716 & 743  & \d_{52} &  6 &   &  4 &\3ol &   & \1ol  &\\
743 & 752 & 716  & \d_{43} &  \6ol &  \5ol &   & &  2 & 1  & \\
\hline
\end{array}
$$
The linear relations $\d_{jk}+\d_{pq}+\d_{rs}=0$ are 
$$
\begin{array}{ccccc}
\d_{23}+\d_{45}+\d_{67}=0  & \qquad & \d_{31}+\d_{46}+\d_{75}=0  & \qquad &\d_{16}+\d_{43}+\d_{52}=0 \\ 
\d_{12}+\d_{65}+\d_{74}=0 & \qquad &   \d_{51}+\d_{37}+\d_{62}=0  & \qquad & \d_{71}+\d_{53}+\d_{24}=0 \\ 
  \phantom{\d_{71}+\d_{53}+\d_{24}=0}  & \qquad & \d_{14}+\d_{27}+\d_{36}=0.  & \qquad & \phantom{\d_{71}+\d_{53}+\d_{24}=0} \\ 
\end{array}
$$

 \begin{prop}
 Let $L$ be an associative 3-plane. Then
 $$
 \{ \d \in \Der_{\sf gr}A_k \; | \; \ker \d = L\} \cong \fso(3,\RR).
 $$
 \end{prop}
\begin{pf}
Using the table we see that 
$$
\{ \d \in \Der_{\sf gr}A_k \;|\; \ker \d =  \langle x_1,x_2,x_3\rangle \} = \langle \d_{12}, \d_{23}, \d_{31} \rangle \cong 
\fso(3,\RR) 
$$
so the result is true for $ \langle x_1,x_2,x_3\rangle$ and hence for $L$ because  $G^c_2$ acts transitively
on associative 3-planes.  
\end{pf}

  \begin{prop}
  \label{prop.deg2.elt}
The element $x_1^2+\cdots+x_7^2$ is $G_2$-invariant and therefore annihilated by $\Der_{\sf gr}A_k$.
\end{prop}
\begin{pf}
 The result follows from the fact that $\{o_1,\ldots,o_7\}$ is an orthonormal basis with respect to the 
 $G_2^c$-invariant inner product $\langle u,v\rangle = \Real(u\vol)$ on $\Im{\mathbb O}$. 
\end{pf}  
 
 \begin{prop} 
Every derivation in $\Der_{\sf gr}A_k$ vanishes on $Q$.
\end{prop}
\begin{pf}
Given the symmetry of $Q$ it suffices to prove that $\d_{23}(Q)=0$. This is a straightforward calculation.
First we use the fact that $\d_{23}$ vanishes on $\langle x_1,x_2,x_3\rangle$. 
Let $i \in \{1,2,3\}$. Because $\d_{23}$ vanishes on $x_i$ and $\d_{23}(x_4)=x_5$ and $\d_{23}(x_5)=-x_4$, $\d_{23}$ vanishes on  
$[x_i,x_4]^2 + [x_i,x_5]^2$. 
Likewise, $\d_{23}$ vanishes on  $[x_i,x_6]^2 + [x_i,x_7]^2$.

Because $\d_{23}(x_4)=x_5$ and $\d_{23}(x_5)=-x_4$, $\d_{23}$ vanishes on $[x_4,x_5]^2$. 
Likewise, $\d_{23}$ vanishes on  $[x_6,x_7]^2$.
It remains to show that $\d_{23}$ vanishes on 
$$
[x_4,x_6]^2 + [x_4,x_7]^2 + [x_5,x_6]^2 + [x_5,x_7]^2.
$$
This is a straightforward calculation. 
\end{pf}

\subsubsection{Remark}
If we identify each $x_i \in A_1$ with $o_i \in \Im {\mathbb O}$, the derivations $\d_{jk}$ become
derivations of ${\mathbb O}$.

\section{The algebra determined by the dual of octonion multiplication}

Let $\nu:{\mathbb O} \otimes {\mathbb O} \longrightarrow {\mathbb O}$ be the multiplication map and 
$\nu^*$ its dual. Define
$$
\wtA:={{T{\mathbb O}^*}\over{(\im \nu^*)}}.
$$
Let $r_1,\ldots,r_7$ be the relations in Proposition \ref{prop.relns1}. 
Then $\wtA$ is the free algebra $k \langle x_0,\ldots,x_7\rangle$ modulo the eight relations
\begin{align*}
\widetilde{r_0}:= & x_0^2-x_1^2-\cdots x_7^2,
\\
\widetilde{r_i}:= & x_ix_0+x_0x_i +r_i, \qquad (1 \le i \le 7). 
\end{align*}
Thus
$$
{{\wtA}\over{(x_0)}} \cong {{A}\over{(x_1^2+\dots +x_7^2)}}.
$$
The algebra $R:=k \langle x_0,\ldots,x_7\rangle/(x_0^2-x_1^2-\cdots x_7^2)$ is one of Zhang's global dimension two, one-relator algebras. In particular, $R$ is coherent and, when $k$ is algebraically closed, 
$\proj R$ is derived equivalent to the Kronecker quiver with eight arrows \cite{M}, \cite{Pion1}.

\end{document}